\documentclass[11pt,a4paper]{article}
\usepackage{authblk}
\usepackage[english]{babel}
\usepackage[T1]{fontenc}
\usepackage[utf8]{inputenc}
\usepackage{float}
\usepackage{amsmath,amsfonts,amsthm,amssymb}
\usepackage{mathtools}
\usepackage{comment}
\usepackage{bbold}
\usepackage{palatino}
\usepackage{eurosym}
\usepackage{enumerate}
\theoremstyle{definition}
\newtheorem{remark}{Remark}
\usepackage[symbol]{footmisc}

\usepackage[toc,page]{appendix}
\usepackage[inline]{enumitem}
\usepackage[%  
    colorlinks=true,
    pdfborder={0 0 0},
    linkcolor=red
]{hyperref}
\usepackage{hyperref}
\hypersetup{colorlinks,linkcolor={blue},citecolor={blue},urlcolor={red}}  
\usepackage[top=2.5cm, left=2.5cm, right=2.5cm, bottom=2.5cm]{geometry}
\usepackage[square,numbers]{natbib}
\numberwithin{equation}{section}
\title{A hybrid discrete-continuum modelling approach to explore the impact of T-cell infiltration on anti-tumour immune response\\
}
\date{\vspace{-6ex}}
\vspace{3ex}
\author[1]{Luis Almeida}
\author[2,3]{Chloe Audebert}
\author[1]{Emma Leschiera\thanks{Corresponding author \\
\textit{Email addresses}: \texttt{luis.almeida@sorbonne-universite.fr} (Luis Almeida), \texttt{chloe.audebert@sorbonne-universite.fr} (Chloe Audebert),  \texttt{emma.leschiera@sorbonne-universite.fr} (Emma Leschiera), \texttt{tommaso.lorenzi@polito.it} (Tommaso Lorenzi).\\
E.L. has received funding from the European Research Council (ERC) under the European Union's Horizon2020 research and innovation programme (grant agreement No 740623). \\
T.L. gratefully acknowledges support from the Italian Ministry of University and Research (MUR) through the grant ``Dipartimenti di Eccellenza 2018-2022'' (Project no. E11G18000350001) and the PRIN 2020 project (No. 2020JLWP23) ``Integrated Mathematical Approaches to Socio–Epidemiological Dynamics'' (CUP: E15F21005420006). \\
L.A., E.L. and T.L. gratefully acknowledge support from the CNRS International Research Project ``Modélisation de la biomécanique cellulaire et tissulaire'' (MOCETIBI).}}
\author[4]{Tommaso Lorenzi}

\affil[1]{\small\textit{Sorbonne Universit\'e, CNRS, Universit\'e de Paris, Inria, Laboratoire Jacques-Louis Lions UMR 7598, 75005 Paris, France.}}
\affil[2]{\textit{Sorbonne Universit\'e, CNRS, Universit\'e de Paris, Laboratoire Jacques-Louis Lions UMR 7598, 75005 Paris, France.}}
\affil[3]{\textit{Sorbonne Universit\'e, CNRS, Institut de biologie Paris-Seine (IBPS), Laboratoire de Biologie Computationnelle et Quantitative UMR 7238, 75005 Paris, France.}}
\affil[4]{\textit{Department of Mathematical Sciences ``G. L. Lagrange'', Dipartimento di Eccellenza 2018-2022, Politecnico di Torino, 10129 Torino, Italy.}}

\begin{document}
%\begin{frontmatter}
\maketitle
\begin{abstract}
We present a spatial hybrid discrete-continuum modelling framework for the interaction dynamics between tumour cells and cytotoxic T cells, which play a pivotal role in the immune response against tumours. In this framework, tumour cells and T cells are modelled as individual agents while chemokines that drive the chemotactic movement of T cells towards the tumour are modelled as a continuum. We formally derive the continuum counterpart of this model, which is given by a coupled system that comprises an integro-differential equation for the density of tumour cells, a partial differential equation for the density of T cells, and a partial differential equation for the concentration of chemokines. We report on computational results of the hybrid model and show that there is an excellent quantitative agreement between them and numerical solutions of the corresponding continuum model. These results shed light on the mechanisms that underlie the emergence of different levels of infiltration of T cells into the tumour and elucidate how T-cell infiltration shapes anti-tumour immune response. Moreover, to present a proof of concept for the idea that, exploiting
the computational efficiency of the continuum model, extensive
numerical
simulations could be carried out, we investigate the impact of T-cell
infiltration on the response of tumour cells to different types of
anti-cancer immunotherapy.
\end{abstract}


\section{Introduction}
\label{introduction}
\subsection{Biological background}
Understanding the cellular processes that underlie the early stages of tumour development and tumour-immune interaction is important to guide the design of effective treatments, especially immunotherapy~\citep{basu2016cytotoxic, boissonnas2007vivo, halle2016vivo, pitt2016targeting}. Experimental and clinical evidence indicates that the immune system plays a critical role in the prevention and eradication of tumours, by detecting immunogenic tumour cells through mutational or abnormally expressed genes and mounting an adaptive immune response~\citep{coulie2014tumour}. In particular, specific immune cells, such as cytotoxic T cells, are activated in secondary lymphoid organs draining the tumour site. Then, these T cells migrate to the tumour micro-environment (TME) in an attempt to eliminate the tumour~\citep{boissonnas2007vivo, miller2003autonomous}. 

However, a myriad of immunosuppressive strategies, the so-called immune checkpoints, help tumour cells  acquiring features that allow them to evade immune detection, which may ultimately result in tumour escape. One important route towards this escape is created as tumour cells highjack the regulatory pathways of the immune system to suppress its functionality~\citep{rabinovich2007immunosuppressive}. The programmed cell death protein-1 (PD1) and its ligand PD-L1 are amongst these inhibitory pathways~\citep{hegde2016and,iwai2002involvement}. Under protracted immune stress, PD-L1 expression can be induced on tumour cells, leading to T cell exhaustion and resistance to anti-tumour immune action in many cancers, such as melanoma~\citep{spranger2015melanoma,wieland2018t} and non-small cell lung cancer (NSCLC)~\citep{mcgranahan2016clonal}. Moreover, the engagement of various oncogenic pathways results in the expression of cytokines and chemokines that mediate the exclusion of T cells from the TME~\citep{kato2017prospects} or, alternatively, the repression of factors that facilitate T cell trafficking and infiltration into the tumour~\citep{spranger2015melanoma}. In this context, the design of immune checkpoint therapies which target regulatory pathways in T cells to enhance anti-tumour immune responses may be beneficial to the treatment of multiple types of cancer~\citep{hegde2016and,kato2017prospects,ribas2018cancer}.

The observation that type, density and location of immune cells within the tumour may be associated with prognosis in different types of cancer led to the development of the ‘immunoscore’ as a prognostic
marker in cancer patients ~\citep{angell2013immune,galon2019approaches, galon2006type, galon2016immunoscore}. The immunoscore provides a score that increases with the density of CD8+ and CD3+ T cells
both in the centre and at the margin of the tumour. CD3 is a protein complex and T cell co-receptor that is involved in activating both cytotoxic T cells (CD8+ naive T cells) and T helper cells (CD4+ naive T cells), and it is therefore the common antigen used to identify both CD4+ and CD8+ T cells. In this vein, a new immune-based, rather than a cancer-based, classification of tumours that relies on
the immunoscore has been proposed in~\citep{galon2006type}, where the authors have classified tumours in four categories.  The ‘‘hot" category comprises tumours which are highly infiltrated by T cells and thus have a high immunoscore. The category ‘‘altered-immunosuppressed" is for tumours with a small amount of infiltrated T cells. Tumours in the ‘‘altered-excluded" category are characterised by two different regions: their margin is T cell-infiltrated while the centre is not. Tumours in these two categories have an intermediate immunoscore. Finally, ‘‘cold" tumours have a low immunoscore and are often associated with a poor response to immunotherapies, since T cells are absent both in the centre of the tumour and at its margin.

\subsection{Mathematical modelling background}
Mathematical models can support a better understanding of the interaction dynamics between tumour cells and immune cells~\citep{eftimie2011interactions,handel2020simulation,makaryan2020modeling,Wilkie2013}. Incorporating in these models the effects of therapeutic strategies that boost anti-tumour immune response can help predicting the success of cancer treatment protocols, including immunotherapy protocols.

Mechanistic dynamical-system models formulated as ordinary differential equations (ODEs) ~\citep{aguade2020tumour,almuallem2021oncolytic,Cattani2010,griffiths2020circulating, konstorum2017addressing, kuznetsov1994nonlinear,Kuznetsov2001,LinErikson2009,luksza2017neoantigen,Takayanagi2001,Wilkie2013} or integro-differential equations (IDEs)~\citep{almeida2021discrete, delitala2013recognition, Kolev2003, lorenzi2015mathematical} have been developed to investigate the interaction dynamics between tumour and immune cells in various scenarios, including immunotherapy. However, these models often rely on the assumption that
cells are well-mixed and, therefore, do not take into account spatial dynamics of tumour and immune cells. As a result, partial differential equation (PDE) models to study the spatio-temporal dynamics of tumour-immune interactions have also been developed -- see for example~\citep{al2012evasion,atsou2020size}. 

Deterministic continuum models are amenable to both analytical and numerical approaches, which allow for a more in-depth theoretical understanding of the underlying cellular dynamics. However, they make it more difficult to incorporate the finer details of dynamical interactions between tumour cells and T cells. Moreover, they cannot easily capture the emergence of population-level phenomena that are driven by stochastic fluctuations in single-cell biophysical properties, which are relevant in the regime of low cellular densities~\citep{lorenzi2022cancer}. Hence, one ideally wants to derive these models as the appropriate limit of stochastic discrete (i.e. individual-based) models of the interaction dynamics between tumour cells and T cells. A number of stochastic discrete models~\citep{christophe2015biased,macfarlane2018modelling,macfarlane2019stochastic} and hybrid discrete-continuum models~\citep{kim2012modeling,mallet2006cellular,leschiera2022mathematical} have also been used to study the interaction dynamics between tumour and immune cells. In contrast to continuum models, these
discrete models track the dynamics of single cells, thus permitting the representation of single cell-scale mechanisms, and account for possible stochastic fluctuations in single-cell biophysical properties. Integrating the results of computational simulations of stochastic discrete models with the results of analysis and numerical simulation of their deterministic continuum counterparts makes it possible to identify more clearly the validity domain of the results obtained, thus leading to more robust biological insights. As a consequence, the derivation of continuum models for the dynamics of cell populations from underlying hybrid models has become an active research field -- see, for instance, ~\citep{byrne2009individual,champagnat2008individual,chisholm2016evolutionary,johnston2015modelling}. 

\subsection{Contents of the article}
%The rules governing cell dynamics in the hybrid model result in a discrete-time branching random walk on a regular lattice~\citep{hughes1995random}. 
In this article, we develop a spatial hybrid discrete-continuum model for the interaction dynamics between tumour cells and immune cells. In this framework, a stochastic individual-based model tracking the dynamics of single tumour cells and immune cells is coupled with a balance equation for the concentration of chemokines (\textit{e.g.} CXCL9/10) which are secreted by tumour cells and drive the chemotactic movement of immune cells towards the tumour~\citep{galon2019approaches,gorbachev2007cxc}. While being aware of the fact that a variety of different types of cells and molecules take part in the immune response against tumours, here we focus on cytotoxic T cells only, since they are the immune agents that are most commonly stimulated by immunotherapies~\citep{huang2017t,tumeh2014pd}. 

%In particular, possible alterations in T-cell motility observed at high cell densities (i.e. volume-filling effects)~\citep{painter2002volume} are taken into account by modulating the probabilities of T-cell movement by a decaying function of the densities of tumour cells and T cells. 
In this model, cell dynamics are governed by a set of rules that result in a discrete-time branching random walk on a regular lattice~\citep{hughes1995random}. Using methods similar to those we have previously employed in~\citep{almeida2021discrete,bubba2020discrete,macfarlane2020hybrid}, we formally derive the deterministic continuum counterpart of the hybrid model, which is given by a coupled system that comprises an IDE for the density of tumour cells, a PDE for the density of T cells, and a PDE for the concentration of chemokines. We report on computational results of the hybrid model and show that there is an excellent quantitative agreement between them and numerical solutions of the corresponding continuum model. These results shed light on the mechanisms that underlie the emergence of different levels of infiltration of T cells into the tumour and elucidate how T-cell infiltration shapes anti-tumour immune response. Moreover, to present a proof of concept for the idea that, exploiting the computational efficiency of the continuum model, extensive numerical simulations could be carried out to identify possible targets to improve the efficacy of anti-cancer therapy, we investigate the impact of T-cell infiltration on the dynamics of tumour cells under parameter settings which provide a simplified representation of the action of different types of  immunotherapy.

The article is organised as follows. In Section~\ref{sec:hybrid model}, the hybrid discrete-continuum model is introduced. In Section~\ref{sec:Corresponding continuum model},
the deterministic continuum counterpart of this model is presented (a formal derivation is provided
in Appendix~\ref{Formal derivation of the continuum model corresponding to the hybrid model}). In Section~\ref{sec:Numerical simulations}, computational results of the
hybrid model are discussed and integrated with numerical solutions of the continuum model. In Section~\ref{sec:Conclusive discussion and research perspectives}, biological implications of the main findings of this study are summarised and directions
for future research are outlined.

%We assume that T cells in the TME have a chemotactic behaviour, \textit{i.e.} they move towards zones of high concentrations of the chemoattractant, which is secreted by tumour cells and diffuses in the spatial domain. The chemotactic movement of T cells plays a fundamental role in the anti-tumour immune response; it directs T cell trafficking and infiltration into the tumour site and it is crucial to understand some of the mechanisms that lead to immune evasion~\citep{galon2019approaches}. Volume-filling effects, which consider volume limitations and the finite cell size on the ease of movement \citep{painter2002volume}, are taken into account by modulating the probabilities of T cell movement by a decaying function of the densities of tumour cells and T cells. This captures possible alterations in the spatial dynamics leading to the infiltration of T cells into the tumour. 

%, which is coupled with a discrete balance equation for the concentration of the chemoattractant. The modelling approach adopted here to describe T cell movement is inspired from , where also volume-filling effects with cell population interaction
% have been considered. 

\section{Hybrid discrete-continuum model}
\label{sec:hybrid model}
In our model, each cell is seen as an agent that occupies a position on a lattice, while the concentration of chemokines, to which we will refer as ``chemoattractant'' in the remainder of the article, is described by a discrete, non-negative function. Each tumour cell can proliferate or die at certain rates. In the vein of~\citep{atsou2020size, leschiera2022mathematical,macfarlane2019stochastic}, here we focus on tumours in the early stages of development (\textit{i.e.} small pre-angiogenic tumours) and, therefore, we neglect the effects of the movement of tumour cells.  As similarly done in~\citep{atsou2020size}, we let the chemoattractant be secreted by tumour cells and undergo both natural decay and linear diffusion. T cells enter the spatial domain where the tumour is located through blood vessels at a rate that is proportional to the total amount of chemoattractant. Upon entering the domain, T cells undergo undirected, random movement and chemotactic movement towards regions of higher concentration of the chemoattractant (\textit{i.e.} cells migrate towards the tumour), and exert a cytotoxic action against tumour cells.

%but there would be no additional difficulty in considering a three-dimensional domain
For ease of presentation, we let the cells and the chemoattractant be distributed across a $d$-dimensional domain $\Omega$, with $d=1$ or $d=2$. In particular, we consider the case where the spatial domain is represented by the set $\Omega:=[0, \ell]^d$, with $\ell\in \mathbb{R^*_+}$, where $\mathbb{R^*_+}$ is the set of positive real numbers not including zero. The position of the cells and the molecules of chemoattractant at time $t \in [0,t_f] \subset \mathbb{R}_+$ is modelled by the variable $x \in [0, \ell]$ when $d=1$ and by the vector $\mathbf{x}=(x,y)\in [0, \ell]^2$ when $d=2$. %Let $[x]$ denote the unit measuring the position occupied by cells. \\

We discretise the time variable $t$ and the space variables $x$ and $y$, respectively, as $t_k =k \tau$, $x_{i} =i\chi$ and $y_{j} =j\chi$, with $k \in \mathbb{N}_0$ and $(i,j) \in [0, \mathcal{N}]^2 \subset \mathbb{N}_0^2$, where $\mathbb{N}_0$ is the set of natural numbers including zero. Here, $\tau \in \mathbb{R^*_+}$ and $\chi \in \mathbb{R^*_+}$ are the time- and space-step, respectively, and $\mathcal{N} := 1 + \lceil\frac{\ell}{\chi}\rceil$, where
$\lceil\cdot \rceil$ denotes the ceiling function. We will use the notation $\mathbf{i}\equiv i$ and $\mathbf{x_i}\equiv x_i$ when $d=1$, and $\mathbf{i}\equiv (i,j)$ and $\mathbf{x_i}\equiv (x_i,y_j)$ when $d=2$.\\

We denote by $n_\mathbf{i}^k$ the density of tumour cells, which is defined as the number of tumour cells at position $\mathbf{x_i}$ and at time $t_k$, $N_\mathbf{i}^k \in \mathbb{N}_0$, divided by the size of a lattice site, that is
\begin{equation}
    n_\mathbf{i}^k\equiv n(\mathbf{x_i},t_k):=\frac{N_\mathbf{i}^k}{\chi^d}.
    \label{eq:tumourDensity}
\end{equation}
%Let $[cell]$ denotes the number of cells. The density of tumour cells $n$ is then measured in $[cell]\cdot[x]^{-1}$. \\

Furthermore, we denote by $c_\mathbf{i}^k$ the density of T cells, which is defined as the number of T cells at position $\mathbf{x_i}$ and at time $t_k$, $C_\mathbf{i}^k \in \mathbb{N}_0$, divided by the size of a lattice site, that is
\begin{equation}
    c_\mathbf{i}^k\equiv c(\mathbf{x_i},t_k):=\frac{C_\mathbf{i}^k}{\chi^d}.
     \label{eq:tcellDensity}
\end{equation}
%As the tumour cell density, the density of T cells $c$ is  measured in $[cell]\cdot[x]^{-1}$. \\

Finally, the concentration of chemoattractant on the lattice site $\mathbf{i}$ and at time-step $k$ is modelled by the discrete, non-negative function $\phi_\mathbf{i}^k\equiv \phi(\mathbf{x_i},t_k)$. %Having in mind the chemical nature of the signal, the attractive potential can be measured in number of chemoattractant molecules, with a unit denoted by $[mol]$. Therefore, the concentration of the chemoattractant $\phi_\mathbf{i}^k$ is measured in $[mol]\cdot[x]^{-1}.$ Notice that the dynamics of the chemoattractant is modelled deterministically, rather than as a random walk. 

\begin{figure}
      \centering
      \includegraphics[height=8cm]{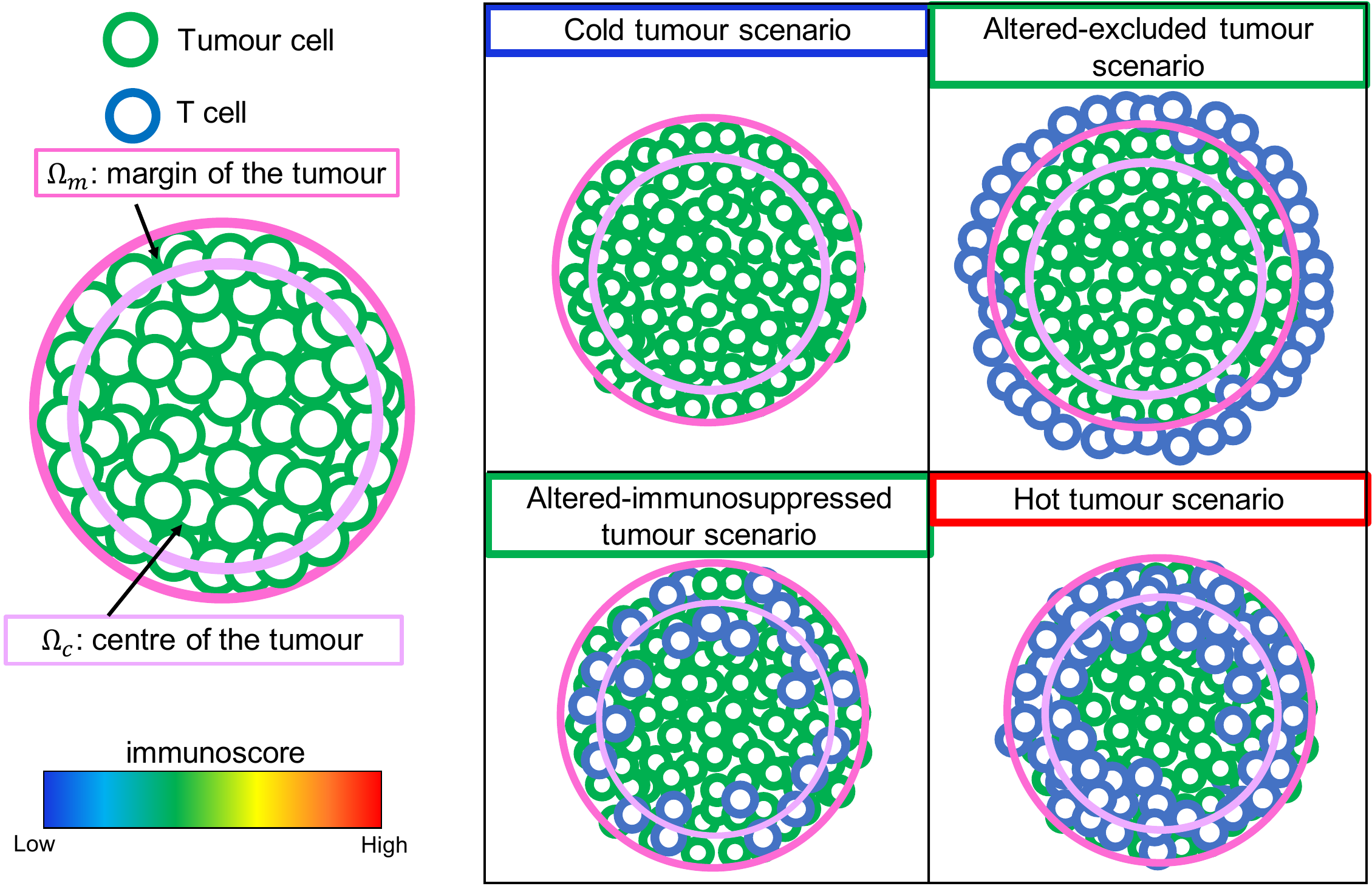}
      \caption{\textbf{Schematics describing the classification of cold,  altered-excluded, altered-immunosuppressed and hot tumour scenarios based on the immunoscore.} In the mathematical framework of our model we classify different tumour scenarios
depending on the value of the immunoscore at the end of numerical simulations, \textit{i.e.} the quantity defined via \eqref{eq:immunoscore_fin}. The immunoscore is here defined as a weighted sum of the number of T cells in the centre
of the tumour and the number of T cells at the margin
of the tumour. In more detail, along the lines of \cite{galon2019approaches}, scenarios for which the value of the immunoscore is low will be classified as ‘cold
tumour scenarios’; scenarios with an intermediate value of the immunoscore will be classified
as ‘altered tumour scenarios’, which will then be further classified as ‘altered-immunosuppressed tumour scenarios’ or ‘altered-excluded tumour scenarios’ based
on the distribution of T cells in the centre and at the margin of the tumour; finally, scenarios characterised by a high value of the immunoscore will be classified as ‘hot tumour scenarios’.}
      \label{fig:immunoscore_schema}
  \end{figure}

In the mathematical framework of our model, the quantity
\begin{equation}
    I^k\equiv I(t_k):=\frac{\lvert\Omega_\text{c}\rvert}{\lvert\Omega_{\text{tum}}\lvert}I_{\Omega_\text{c}}(t_k)+\frac{\lvert\Omega_\text{m}\lvert}{\lvert\Omega_{\text{tum}}\lvert}I_{\Omega_\text{m}}(t_k)
    \label{eq:immunoscore}
\end{equation}
provides a possible simplified measure of the immunoscore $I$ at time $t_k$. In~\eqref{eq:immunoscore}, $I_{\Omega_\text{c}} \in \mathbb{N}_0$ is the number of T cells within the set $\Omega_\text{c}\subset\Omega$ defined as the ‘centre of the tumour’, $I_{\Omega_\text{m}} \in \mathbb{N}_0$ is the number of T cells within the set $\Omega_\text{m}\subset \Omega$ defined as the ‘margin of the tumour’, and $\Omega_{\text{tum}} := \Omega_\text{c}\,\cup\,\Omega_\text{m}$ is the whole region occupied by the tumour. Here, $\lvert(\cdot)\rvert$ is the measure of the set $(\cdot)$. %, which corresponds to the size of the interval $(\cdot)$ if $d=1$ or the area of the square $(\cdot)^2$ if $d=2$. 
Given the initial distribution of tumour cells, we define ({\it cf.} Fig.~\ref{fig:immunoscore_schema}) %To give more relevance to the T cells in the centre of the tumour than those at the margin of it, in Equation \eqref{eq:immunoscore} we choose $\Omega_c$ and $\Omega_m$ such that $|\Omega_\text{m}|<|\Omega_\text{c}|$.
\begin{equation}
    \Omega_\text{c}:=\left\{\mathbf{x_i}\in \Omega : \lVert\mathbf{x_i}-\mathbf{x}_{\text{cm}}\rVert_d<R \right\} \quad \text{and} \quad I_{\Omega_\text{c}}(t_k):=\sum_{\mathbf{i}} c_\mathbf{i}^k\,\chi^d\,\mathbb{1}_{\Omega_\text{c}}(\mathbf{x_i})
    \label{eq:centreTum}
    \end{equation}
while
\begin{equation}
    \Omega_\text{m}:=\Omega_{\text{tum}}\setminus \Omega_\text{c} \quad \text{and} \quad I_{\Omega_\text{m}}(t_k):=\sum_{\mathbf{i}}c_\mathbf{i}^k\,\chi^d\,\mathbb{1}_{\Omega_\text{m}}(\mathbf{x_i}).
    \label{eq:marginTum}
    \end{equation}
    Here, $\mathbf{x}_{\text{cm}}\in \Omega$ is the initial centre of mass of the tumour, which is computed  as 
    \begin{equation}
    \mathbf{x}_{\text{cm}}=\frac{1}{\rho_n^0}\sum_{\mathbf{i}} n_\mathbf{i}^0 \chi^d \mathbf{x_i},
    \end{equation}
    where $\rho_n^0$ is the initial number of tumour cells.
    %$\Omega_c$ represents the interval centered in $\mathbf{x_{cm}}$ and of size $2R$ if $d=1$ or the disk of centre $\mathbf{x_{cm}}$ and radius $R$ if $d=2$. % $||(\cdot)||_d$  
    Moreover, $\mathbb{1}_{\Omega_\text{m}}$ and $\mathbb{1}_{\Omega_\text{c}}$ are the indicator functions of the sets $\Omega_\text{m}$ and $\Omega_\text{c}$, respectively. Note that, in the definition of sets $\Omega_\text{c}$ and $\Omega_\text{m}$,  we are supposing that the radius $R$ is fixed and, therefore, the two sets do not change over time. This is coherent with the fact that, as mentioned earlier, tumour-cell movement is neglected.\\
Abstracting from the ‘immunoscore'-based classification of tumours
recalled in Section~\ref{introduction}, throughout the article we will classify different tumour scenarios depending on the value of $I$ at the end of numerical simulations, \textit{i.e.} the quantity\begin{equation}
    I^f\equiv I(t_f)=\frac{\lvert\Omega_\text{c}\lvert}{\lvert\Omega_{\text{tum}}\lvert}I_{\Omega_\text{c}}(t_f)+\frac{\lvert\Omega_\text{m}\lvert}{\lvert\Omega_{\text{tum}}\lvert}I_{\Omega_\text{m}}(t_f).
   \label{eq:immunoscore_fin}
\end{equation}
In particular, along the lines of \cite{galon2019approaches}, scenarios for which the value of $I^f$ is low will be classified as ‘cold tumour scenarios’; scenarios with an intermediate value of $I^f$ will be classified as ‘altered tumour scenarios’, which will then be further classified as ‘altered-immunosuppressed tumour scenarios’ or ‘altered-excluded tumour scenarios’ based on the distribution of T cells at the centre and margin of the tumour;
finally, scenarios characterised by a high value of $I^f$ will be classified as
‘hot tumour scenarios’. This tumour classification is illustrated by the schematics presented in Fig.~\ref{fig:immunoscore_schema}.
  
\begin{figure}
      \centering
      \includegraphics[height=9cm]{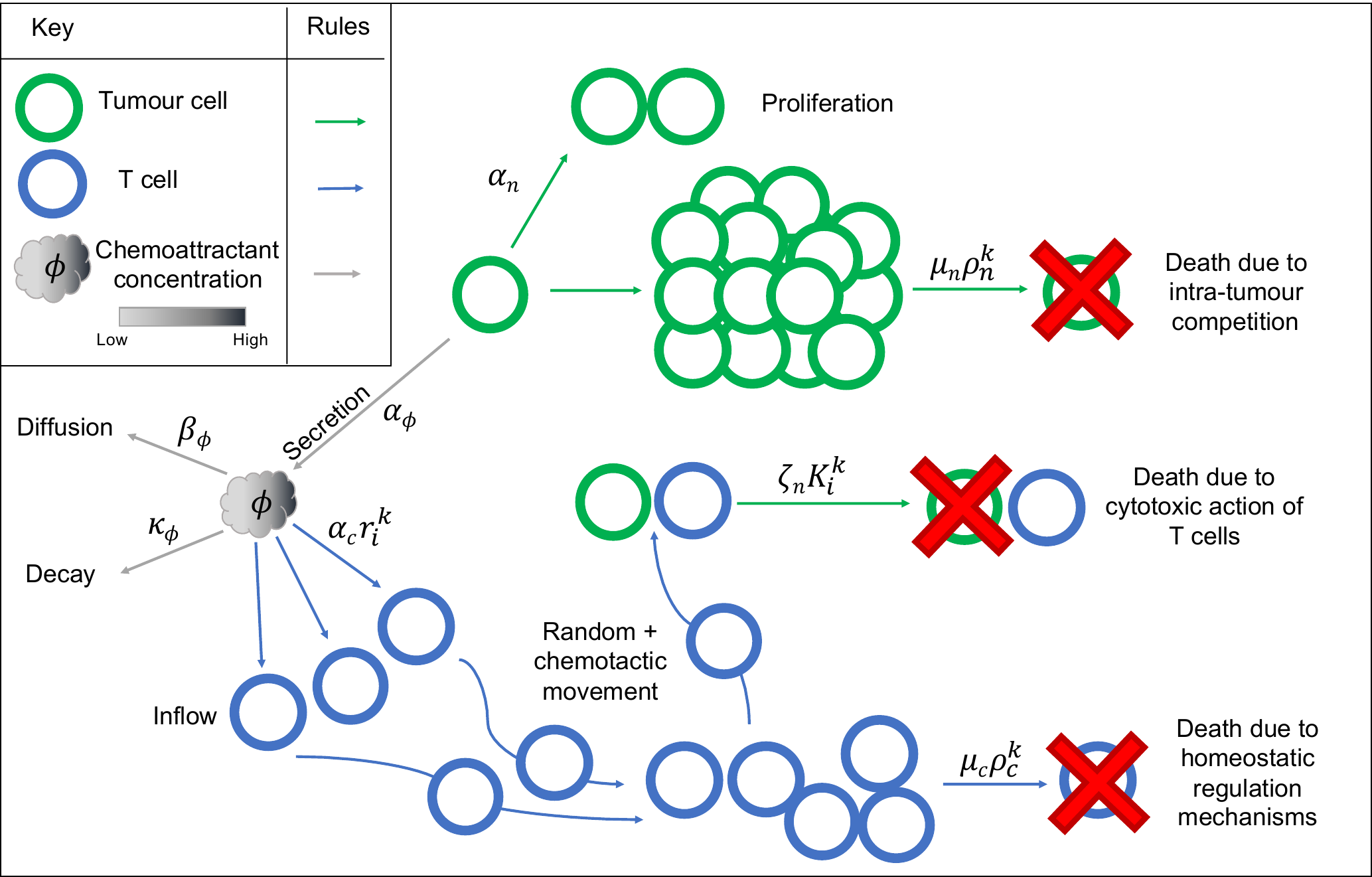}
      \caption{\textbf{Schematics summarising the mechanisms and processes included in the hybrid discrete-continuum model.} The governing rules for the dynamics of tumour cells are here represented by green solid arrows, the rules for the dynamics of T cells are represented by blue solid arrows, and the processes underlying the dynamics of the chemoattractant are represented by grey solid arrows.
At each time-step $t_k$, tumour cells divide at rate $\alpha_n$, die due to intra-tumour competition at rate $\mu_n\rho_n^k$ or die due to the cytotoxic action of T cells at rate $\zeta_nK_i^k$. The quantities $\rho_n^k$ and $K_i^k$ are, respectively, defined via \eqref{eq:rho_n} and \eqref{eq:K_i^k}. Tumour cells also secrete the chemoattractant at rate $\alpha_\phi$. The chemoattractant diffuses through the domain at rate $\beta_\phi$ and decays at rate $\kappa_\phi$. The influx rate of T cells from blood vessels at position $x_i$ is modelled by $\alpha_cr_i^k$, with $r_i^k$ defined via \eqref{eq:inflow}. T cells change their position through a combination of undirected, random movement and chemotactic movement in response to the chemoattractant secreted by tumour cells, while tumour cell movement is neglected. Finally, T cells die due to homeostatic regulation mechanisms at rate $\mu_c\rho_c^k$, with $\rho_c^k$ defined via~\eqref{rho_c}. }
      \label{fig:my_label}
  \end{figure}

The strategies used to model the dynamics of the cells and the chemoattractant when $d=1$ are described in detail in the following subsections, and are also schematically illustrated in Fig.~\ref{fig:my_label}. Analogous strategies are used in the case where $d=2$.
 
 In the remainder of the article, when necessary, the subscripts $n$, $c$ and $\phi$ will be used to identify the parameters and functions related to the dynamics of tumour cells, T cells and chemoattractant, respectively. 
  
\subsection{Dynamics of tumour cells}
We consider a scenario where tumour cells proliferate (\emph{i.e.} undergo cell division) and die due to intra-tumour competition as well as to the cytotoxic action of T cells. We assume that a dividing tumour cell is replaced by two identical cells that are placed on the same lattice site as their parent, while a dying cell is removed from the system. \\
\subsubsection{Tumour cell proliferation and death induced by intra-tumour competition}
\label{Tumour cell proliferation and natural death}
At every time-step $k$, we allow tumour cells to undergo cell division with probability 
\begin{equation}
  \tau\alpha_n>0,  
  \label{eq:tumProl}
\end{equation} where $\alpha_n\in\mathbb{R^*_+}$ %measured in $[t]^{-1}$, 
represents the rate of tumour cell proliferation.

In order to capture the effect of cell death induced by intra-tumour competition, we let tumour cells die at a rate proportional to their number, which is denoted by  
\begin{equation}
\label{eq:rho_n}
    \rho_n^k\equiv \rho_n(t_k):=\sum_{i}n_i^k\,\chi,\end{equation}
    with constant of proportionality $\mu_n\in \mathbb{R^*_+}$. Hence, between the time-step $k$ and the time-step $k+ 1$, we let a tumour cell die due to intra-tumour competition with probability  
\begin{equation}
    \tau \,\mu_n \rho_n^k\geq 0.
    \label{eq:ProbmortTum}
\end{equation}

\subsubsection{Cytotoxic action of T cells against tumour cells}
\label{Effect of the immune system on the tumour}
When T cells are sufficiently close to a tumour cell, they release cytotoxic substances which may lead to the death of the tumour cell~\citep{boissonnas2007vivo}. Therefore, building on the modelling strategies proposed in~\citep{atsou2020size}, we let a tumour cell die due to the cytotoxic action of T cells at a rate proportional to the number of T cells in a sufficiently close neighbourhood of the tumour cell. This models the fact that T cells can interact with a tumour cell up to a certain distance, and that beyond such a distance the tumour cell can no longer be induced to death. Hence, we introduce $K_i^k\in \mathbb{N}_0$ which represents, at each time-step $t_k$, the number of T cells that can exert a cytotoxic action against a tumour cell at position $x_i$. In particular, we define $K_i^k$ as
	\begin{equation} 
	K_i^k\equiv K(x_i,t_k):=\sum_{p} \eta(x_i,x_p;\theta)c_p^k\,\chi.
	\label{eq:K_i^k}
	\end{equation}
	The function $\eta$ is defined as follow
 
	    		\begin{equation}
	\eta(x,z;\theta):=\begin{cases}1 \qquad \;\; \text{if} \qquad \vert x-z\vert\leq\theta \\
	    	0  \qquad  \, \; \, \text{if} \qquad \vert x-z\vert>\theta 
	    	\end{cases}
	    	\quad
\text{for} \; (x, z; \theta) \in \Omega \times \Omega \times (0, \lvert\Omega\lvert],
	    	\label{eq:function_eta}
	    	\end{equation}
	    	where $\lvert\Omega\lvert$ denotes the size of the set $\Omega$ (\textit{i.e.} $\lvert\Omega\lvert=\ell$ if $\Omega:=[0,\ell]$). The quantity $K_i^k$ defined via \eqref{eq:K_i^k} and \eqref{eq:function_eta} represents the number of T cells within a distance $\theta$ of $x_i$. The parameter $0<\theta\leq \lvert\Omega\lvert$ regulates the maximum radius of interaction between a tumour cell at position $x_i$ and the T cells in its neighbourhood. Therefore, we define the probability of death of tumour cells at position $x_i$ and time $t_k$ due to the cytotoxic action of T cells as 
	  \begin{equation}
     \tau\zeta_n	K_i^k\geq0.
     \label{eq:mortCompetition}
      \end{equation}
      The parameter $\zeta_n\in\mathbb{R^*_+}$ %, measured in $[cell]^{-1}\cdot[t]^{-1}$, 
      is linked to the level of efficiency of T cells at eliminating tumour cells. In particular, lower values of  $\zeta_n$ correspond to scenarios in which this cytotoxic action of T cells is less effective, due for example to a higher expression of PD1 inhibitory receptors and PD-L1 ligands on the surface of T cells and tumour cells~\citep{hegde2016and,iwai2002involvement}.

 \begin{remark}
  Note that \eqref{eq:tumProl}, \eqref{eq:ProbmortTum} and \eqref{eq:mortCompetition} implicitly require the time-step $\tau$ to be sufficiently small that $\tau\alpha_n + \tau \left(\mu_n \rho_n^k + \zeta_n K_i^k \right)$ is less than or equal to 1.    
 \end{remark}

 \subsection{Dynamics of the chemoattractant}
 \label{Dynamic of for the chemoattractant}
We denote by $\phi_i^k$ the concentration of chemoattractant at position $x_i$ and at time $t_k$. The dynamic of $\phi_i^k$ is governed by the following discrete balance equation
\begin{equation}
    \phi_i^{k+1}=\phi_i^k+ \tau \beta_{\phi}(\mathcal{L}\phi^k)_i+\tau\alpha_{\phi}n_i^k-\tau\kappa_{\phi}\phi_i^k, \quad i\in [1,\mathcal{N}-1],
    \label{chemoattr}
\end{equation}
subject to a suitable initial condition and discrete zero-flux boundary conditions, \textit{i.e.} \begin{equation}
 \phi^k_{0}=\phi^k_{1} \; \text{ and } \; \phi^k_{\mathcal{N}}=\phi^k_{\mathcal{N}-1}.
\label{eq:discreteboundaryCondition}
\end{equation}
In the balance equation~\eqref{chemoattr}, $\mathcal{L}$ is the second-order central difference operator on the lattice $\{x_i\}_{i}$, \textit{i.e.}
\begin{equation}
    (\mathcal{L}\phi^k)_i=\frac{1}{\chi^2}\Big(\phi_{i+1}^k+\phi_{i-1}^k-2\phi_{i}^k\Big)
    \end{equation}
Moreover, the parameter $\beta_{\phi}\in \mathbb{R^*_+}$ is the diffusion coefficient of the chemoattractant and % measured in $[x]^2\cdot[t]^{-1}$,
$\kappa_{\phi}\in\mathbb{R^*_+}$ % measured in $[t]^{-1}$, 
is the rate at which the chemoattractant undergoes natural decay. Finally, the parameter $\alpha_{\phi}\in \mathbb{R^*_+}$ represents the per capita production rate of the chemoattractant by tumour cells. The balance equation~\eqref{chemoattr} is simply a standard discretisation of a reaction-diffusion equation of the type that is commonly used to describe the dynamics of molecular species, see for example~\citep{maini1997spatial}. 

\subsection{Dynamics of T  cells}
\label{Dynamics of T  cells}
Following~\citep{gong2017computational}, we consider a scenario where T cells are recruited from different sources corresponding to blood vessels that are located in the tissue surrounding the tumour. 
%We suppose that these sources are equidistant from each other and from the centre of the domain. Moreover, 
T cells can change their position according to a combination of undirected, random movement and chemotactic movement, which are regarded as independent processes. Finally,  T cells can die at a certain rate due to homeostatic regulation mechanisms, and dying cells are removed from the system. This results in the following rules for the dynamics of T cells. 
\subsubsection{Inflow and death of T  cells}
\label{Inflow and death of T  cells}
We let $\omega \subset \Omega$ be the set of points in the tissue surrounding the tumour that are occupied by blood vessels, through which new T cells can enter the domain. Since we do not consider the formation of new blood vessels, we assume that $\omega$ is given and does not change in time.  We denote by $r^k_i$ the term controlling the inflow of T  cells from blood vessels, which is defined as  
\begin{equation}
    r_i^k\equiv r(x_i,t_k):= \phi_{tot}^k\mathbb{1}_\omega(x_i),
    \label{eq:inflow}
\end{equation}
where $\mathbb{1}_\omega$ is the indicator function of the set $\omega$ and $\phi_{tot}^k$ is the total amount of chemoattractant at time $t_k$, that is,
\begin{equation}
    \phi_{tot}^k\equiv \phi_{tot}(t_k):=\sum_{i}\phi_i^k\,\chi.
    \label{eq:phi_tot}
\end{equation}
We then let the influx rate of T cells from blood vessels at position $x_i$ and time $t_k$ be proportional to $r_i^k$ with constant of proportionality $\alpha_c\in \mathbb{R}^*_+$. Hence, between the time-step $k$ and the time-step $k+ 1$, we let a density of T cells equal to 
\begin{equation}
    \tau \alpha_cr^k_i\geq 0
    \label{eq:inflowProb}
\end{equation}
enter a blood vessel at position $x_i$.

Finally, we let T cells die %, measured in $[cell]^{-1}\cdot[t]^{-1}$ 
due to homeostatic regulation mechanisms. In analogy with the case of tumour cells, we suppose the rate of death of T cells to be proportional to the number of T cells  
\begin{equation}
    \rho^k_c\equiv\rho_c(t_k):=\sum_{i}c_i^k\,\chi,
    \label{rho_c}
    \end{equation}
with constant of proportionality $\mu_c\in \mathbb{R}^*_+$.
    Therefore, between the time-step $k$ and the time-step $k+ 1$, we let a T cell die with probability  
\begin{equation}
    \tau \mu_c \rho_c^k\geq 0.
    \label{eq:deathTcells}
\end{equation}

\begin{remark}
Note that \eqref{eq:deathTcells} implicitly requires the time-step $\tau$ to be sufficiently small that the corresponding quantity is less than or equal to 1.
\end{remark}

\subsubsection{Chemotactic movement of T cells}	 
\label{Mathematical modelling of chemotactic T  cell movement}	
We now turn to modelling the chemotactic movement of T cells (\textit{i.e.} the movement of T cells up the gradient of the chemoattractant $\phi_i^k$). Building on~\citep{bubba2020discrete}, chemotactic movement is here modelled as a biased random walk whereby the movement probabilities depend on the difference between the concentration of chemoattractant at the site occupied by a T cell and the concentration of chemoattractant at the neighbouring sites. To take into account possible reduction in cell motility at high cell densities~\citep{painter2019mathematical,slaney2014trafficking,van2017migrating}, we incorporate into the model volume-filling effects~\citep{painter2002volume}, whereby T cell movement is allowed only to site locations $x_i$ where the total cell density $w_i^k:=n_i^k+c_i^k$ is smaller than a threshold value $w_{\text{max}} \in \mathbb{R}^*_+$, which corresponds to a cell tight-packing state. Therefore, we modulate the movement probabilities of T cells by a monotonically decreasing function of the total cell density at the neighbouring sites. Specifically, as similarly done in~\citep{wang2007classical}, we define this function as
\begin{equation}
    \psi(w_i^k):=\begin{cases}
           1-\dfrac{w_i^k}{w_{\text{max}}}, \quad 0\leq w_i^k\leq w_{\text{max}}, \\
           0, \qquad \qquad \quad \text{otherwise}.
    \end{cases}
    \label{sqeezing_prob}
\end{equation}
%We expect the value of the parameter $w_{\text{max}}$ in definition~\eqref{sqeezing_prob} to affect the level of T-cell infiltration into the tumour. In fact, large values of $w_{\text{max}}$ enable T-cell movement even at high cell densities, thus allowing higher levels of T-cell infiltration to be attained. On the contrary, low values of
%$w_{\text{max}}$ lead T-cell movement to stop even at low cell densities, which may result in the accumulation of T cells at the margin of the
%tumour and lower levels of T-cell infiltration.\\
Hence, 
for a T cell on the lattice site $x_i$ and  at time step $t_k$, we define: 
\begin{enumerate}
    \item the probability of moving to the lattice site $x_{i-1}$ (\textit{i.e.} the probability of moving left) via chemotaxis as		
\begin{equation}
    q_{Li}^k:=\nu\,\psi(w_{i-1}^k)\Big(\frac{\phi_{i-1}^k-\phi_i^k}{2\phi_{\max}}\Big)_{+} \quad \text{for} \quad i \in [1,\mathcal{N}] \quad \text{and} \quad q_{L0}^k=0
    \label{chem_left}
\end{equation}
where  $(\cdot)_+$ denotes the positive part of $(\cdot)$ and $\phi_{\max} \in \mathbb{R}_+^*$ is the maximum value which can be attained by the concentration of the chemoattractant -- see comments below definition~\eqref{chem}; 

\item the probability of moving to the lattice site $x_{i+ 1}$ (\textit{i.e.} the probability of moving right) via chemotaxis as
\begin{equation}
    q_{Ri}^k:=\nu\,\psi(w_{i+1}^k)\Big(\frac{\phi_{i+1}^k-\phi_i^k}{2\phi_{\max}}\Big)_{+} \quad \text{for} \quad i \in [0,\mathcal{N}-1] \quad \text{and} \quad q_{R\mathcal{N}}^k=0;
     \label{chem_right}
\end{equation}
 \item and the probability of not undergoing chemotactic movement as
\begin{equation}
    1-q_{Li}^k-q_{Ri}^k \quad \text{for} \quad i \in [0,\mathcal{N}].
     \label{chem}
\end{equation}
\end{enumerate}
Here, the parameter $\nu\in \mathbb{R_+}$, with $0\leq \nu\leq 1$, is directly proportional to the chemotactic sensitivity of T cells.  Dividing by $\phi_{\max}$ ensures that the values of the quotients in \eqref{chem_left} and \eqref{chem_right} are less than or equal to $1$.

\subsubsection{Undirected, random movement of T cells}
\label{Mathematical modelling of undirected, random T  cell movement}
	To model the effect of undirected, random movement, we allow T  cells to update their position according to a random walk with movement probability $\lambda\in \mathbb{R_+^*}$, where $0<\lambda\leq 1$. In particular, we assume that a T cell on the lattice site $x_i$ can move via undirected, random movement into either the lattice site $x_{i-1}$ or the lattice site $x_{i+1}$ with probability $\dfrac{\lambda}{2}$. As similarly done in the case of chemotactic movement, in order to capture a possible reduction in T cell motility at higher cell densities~\citep{painter2019mathematical,slaney2014trafficking,van2017migrating}, we modulate the movement probability by a decreasing function of the density of tumour cells and T cells at the neighbouring sites. In particular, for a T cell on the lattice site $i$ and at the time-step $k$, we define:
	\begin{enumerate}
	    \item 
	the probability of moving to the lattice site $i-1$ via undirected, random movement as 
	\begin{equation}
   T _{Li}^k:=\frac{\lambda}{2}\,\psi(w_{i-1}^k) \quad \text{for} \quad i \in [1,\mathcal{N}] \quad \text{and} \quad T_{L0}^k=0;
    \label{diff_left}
\end{equation}
\item the probability of moving to the lattice site $i+1$ via undirected, random movement as 
	\begin{equation}
   T _{Ri}^k:=\frac{\lambda}{2}\,\psi(w_{i+1}^k) \quad \text{for} \quad i \in [0,\mathcal{N}-1] \quad \text{and} \quad T_{R\mathcal{N}}^k=0;
    \label{diff_right}
\end{equation}
\item and the probability of not undergoing undirected, random movement as 
\begin{equation}
    1-T_{Li}^k-T_{Ri}^k \quad \text{for} \quad i \in [0,\mathcal{N}].
     \label{diff}
\end{equation}
\end{enumerate}
In \eqref{diff_left} and \eqref{diff_right}, the modulating function $\psi$ is defined via  \eqref{sqeezing_prob}.
\section{Corresponding continuum model}
\label{sec:Corresponding continuum model}
Letting the time-step $\tau \rightarrow 0$ and the space-step
$\chi \rightarrow 0$ in such a way that
\begin{equation}
    \dfrac{\lambda\chi^2}{2d\,\tau}\rightarrow \beta_c \in \mathbb{R}_+^* \quad \text{and} \quad  \dfrac{\nu\chi^2}{2d\,\phi_{\max}\,\tau}\rightarrow \gamma_c \in \mathbb{R}_+^* \quad \text{as } \tau \rightarrow 0, \, \chi \rightarrow 0,
    \label{eq:conditionsContMod}
\end{equation}
using
the method employed in~\citep{almeida2021discrete, bubba2020discrete, macfarlane2020hybrid}, it is possible to formally show (see Appendix~\ref{Formal derivation of the continuum model corresponding to the hybrid model}) that the deterministic continuum counterpart of the hybrid
model described in Section~\ref{sec:hybrid model} comprises the following coupled IDE-PDE-PDE system for the density of tumour cells, $n(\mathbf{x},t)$, the density of T cells, $c(\mathbf{x},t)$ and the chemoattractant concentration, $\phi(\mathbf{x},t)$
\begin{equation}
    \begin{cases}
     \partial_t n =\alpha_n n-\mu_n\rho_n(t)\, n-\zeta_nK(\mathbf{x},t) n %\quad (x,t)\in \left[0,l\right]^d\times\mathbb{R}^*_+
    \\\\
     \partial_t c- \nabla
     \cdot\Big[\beta_c \, \psi(w) \, \nabla c-\gamma_c \, \psi(w) \, c \, \nabla \phi -\beta_c \, c \, \psi^\prime(w) \, \nabla w\Big]=  -\mu_c\rho_c(t)c+\alpha_cr(\mathbf{x},t) %\quad (x,t)\in (0,l)^d\times\mathbb{R}^*_+
     \\\\
    \partial_t \phi-\beta_{\phi}\Delta \phi=\alpha_{\phi}n-\kappa_{\phi}\phi %\quad (x,t)\in (0,l)^d\times\mathbb{R}^*_+\\
    \\\\
    w(\mathbf{x},t):=n(\mathbf{x},t)+c(\mathbf{x},t),
    \end{cases}
    \label{continuum_model}
\end{equation}
where the IDE \eqref{continuum_model}$_1$ is posed on $\Omega\times(0,t_f]$, while the PDEs~\eqref{continuum_model}$_2$ and \eqref{continuum_model}$_3$ are posed on $\Omega\setminus\partial \Omega\times(0,t_f]$ and are subject to zero-flux boundary conditions on $\partial \Omega$. The IDE-PDE-PDE system \eqref{continuum_model} is complemented with the following definitions 
$$K(\mathbf{x},t):=\displaystyle\int_\Omega \eta(\mathbf{x},\mathbf{x^\prime};\theta)c(\mathbf{x^\prime},t)\,\mathrm{d}\mathbf{x^\prime},\quad r(\mathbf{x},t):=\phi_{tot}(t)\mathbb{1}_\omega(\mathbf{x})$$
$$\rho_n(t):=\displaystyle\int_\Omega n(\mathbf{x},t)\,\mathrm{d}\mathbf{x}, \quad \rho_c(t):=\displaystyle\int_\Omega c(\mathbf{x},t)\,\mathrm{d}\mathbf{x}, \quad \phi_{tot}(t):=\displaystyle\int_\Omega \phi(\mathbf{x},t)\,\mathrm{d}\mathbf{x}.$$
In system \eqref{continuum_model}, $\beta_c\in \mathbb{R}_+^*$ defined via \eqref{eq:conditionsContMod} is the diffusion coefficient (\textit{i.e.} the motility) of T cells, while $\gamma_c \in \mathbb{R}_+^*$ defined via \eqref{eq:conditionsContMod} is the chemotactic sensitivity of T cells to the chemoattractant. %Notice that as $\lambda$ and $\eta$ are dimensionless parameters, $\beta_c$ and $\gamma_c$ are respectively measured in $[x]^2\cdot[t]^{-1}$ and $ [x]^2\cdot[t]^{-1}\cdot[mol]^{-1}$. 
\begin{comment}
Similarly, in the case where the cells and the chemoattractant are distributed over a two-dimensional domain $\Omega^2$ of step $\chi$, defining the cell densities via the two-dimensional analogue of \eqref{eq:tumourDensity} and \eqref{eq:tcellDensity}, letting the operator $\mathcal{L}$ be the finite-difference Laplacian on the grid $\{x_{h}\}\times\{y_{i}\}$, $h ,i\in[0,I]$, and assuming $\tau \rightarrow 0$ and
$\chi \rightarrow 0$ in such a way that
\begin{equation}
    \frac{\lambda}{4}\frac{\chi^2}{\tau}\rightarrow \beta_c \in \mathbb{R}_*^+ \quad \text{and} \quad  \frac{\eta}{4\phi_{\max}}\frac{\chi^2}{\tau}\rightarrow \gamma_c \in \mathbb{R}_*^+,
    \label{eq:conditionsContMod2D}
\end{equation}
it is possible to formally obtain the 2D continuum counterpart of our hybrid model, posed on $\Omega^2\times \mathbb{R}^*_+$.
\end{comment}
\section{Numerical simulations}
\label{sec:Numerical simulations}
In this section, we report on computational results of the hybrid discrete-continuum model along with numerical solutions
of the corresponding continuum model given by the IDE-PDE-PDE system \eqref{continuum_model}. First, we establish a baseline scenario in which the level of efficiency of T cells at eliminating tumour cells (\emph{i.e.} the parameter $\zeta_n$) is sufficiently high so as to lead to tumour eradication. Then, we reduce the level of T cell efficiency in order to avoid tumour eradication, and we explore the mechanisms that underlie the emergence of different levels of infiltration of T cells into the tumour, which correspond to cold, altered-immunosuppressed, altered-excluded and hot tumour scenarios. In particular, we carry out sensitivity analysis to two parameters that we expect to play a key role in determining the spatial distribution of T cells: the secretion rate of chemoattractant by tumour cells (\emph{i.e.} the parameter $\alpha_{\phi}$) and the threshold value of the total cell density above which T cell movement is impaired (\emph{i.e.} the parameter $w_{\max}$ in definition~\eqref{sqeezing_prob}). Finally, exploiting the computational efficiency of the continuum model, we investigate the impact of T-cell infiltration on the dynamics of tumour cells under parameter settings which provide a simplified representation of the action of different types of  immunotherapy.

\subsection{Set-up of numerical simulations}
\begin{figure}
\centering
\includegraphics[scale=0.2]{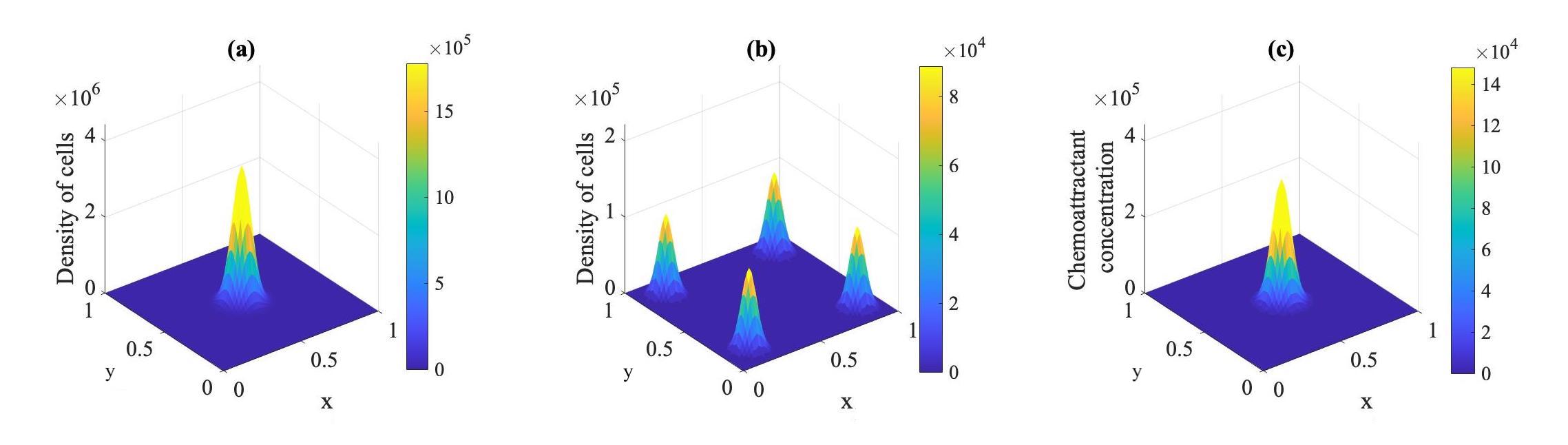}
\caption{\textbf{Initial conditions.} Plots of the density of tumour cells $n(x,y,t)$ (panel \textbf{(a)}), the density of T cells $c(x,y,t)$ (panel \textbf{(b)}), and the chemoattractant concentration $\phi(x,y,t)$ (panel \textbf{(c)}) at the initial time of the simulations (\textit{i.e.} at $t=0$) for the continuum model. Analogous initial conditions are used for the hybrid model.}
\label{fig:initCond}
\end{figure}

The hybrid and continuum models are parameterised using parameter values retrieved from the literature, wherever possible. The full list of parameter values and related references are provided in Table~\ref{table2}.
For the numerical simulations we report on, we use the 2D spatial domain $\Omega:=[0,1]^2$. Under the parameter choice of Table~\ref{table2}, this is
equivalent to considering a square region of a 2D cross-section of a tumour tissue of area $1 \,cm^2$. Furthermore, to carry out numerical simulations of the hybrid model, we use the space-step $\chi= 0.016$~{\it cm} and the time-step $\tau = 1 \times 10^{-4}$~{\it days}. Finally, unless otherwise specified, we choose the final time $t_f = 15$ {\it days}. All simulations are performed in \textsc{Matlab}~\citep{MATLAB:2021}.
\paragraph{Initial conditions and blood vessel distribution}
Fig.~\ref{fig:initCond} displays the initial conditions chosen to carry out numerical simulations.
In particular, for the hybrid model we have 
\begin{equation}
n_{ij}^0=800\exp[-200(i\chi-x_{11}^*)^2-200(j\chi-y_{11}^*)^2],
\label{eq:init_tum}
\end{equation}
\begin{equation}
c_{ij}^0=60\sum_{p=1}^4\exp[-300(i\chi-x_{2p}^*)^2-300(j\chi-y_{2p}^*)^2],
\label{eq:init_imm}
\end{equation}
and
\begin{equation}
\phi_{ij}^0=90\exp[-200(i\chi-x_{11}^*)^2-200(j\chi-y_{11}^*)^2],
\label{eq:init_chem}
\end{equation}
with $(x_{11}^*,y_{11}^*)=(0.5,0.5)$, $(x_{21}^*,y_{21}^*)=(0.26,0.74)$,  $(x_{22}^*,y_{22}^*)=(0.26,0.26)$,  $(x_{23}^*,y_{23}^*)=(0.74,0.74)$ and  $(x_{24}^*,y_{24}^*)=(0.74,0.26)$. The points where T cells are initially concentrated are the centres of the regions where blood vessels are assumed to be located, that is, the set $\omega$ in~\eqref{eq:inflow} is defined as the union of four balls of small radius and centres $(x_{21}^*,y_{21}^*)$, $(x_{22}^*,y_{22}^*)$, $(x_{23}^*,y_{23}^*)$ and  $(x_{24}^*,y_{24}^*)$.  

Similarly, for the continuum model we have
\begin{equation}
n(x,y,0)=800\exp[-200(x-x_{11}^*)^2-200(y-y_{11}^*)^2],
\label{eq:init_tum_cont}
\end{equation}
\begin{equation}
c(x,y,0)=60\sum_{p=1}^4\exp[-300(x-x_{2p}^*)^2-300(y-y_{2p}^*)^2], 
\label{eq:init_imm_cont}
\end{equation}
and
\begin{equation}
\phi(x,y,0)=90\exp[-200(x-x_{11}^*)^2-200(y-y_{11}^*)^2].
\label{eq:init_chem_cont}
\end{equation}
A description of the algorithmic rules that underlie computational simulations of the hybrid model is provided in Appendix~\ref{app:sim1Ddiscrete}, while the methods employed to numerically solve the IDE-PDE-PDE system \eqref{continuum_model}, subject to suitable initial conditions and no-flux boundary conditions, are detailed in Appendix~\ref{app:sim1Dcontinuum}.\\

Given the initial conditions of the two models, we compute the coordinates of the centre of mass of the tumour $(x_\text{cm},y_\text{cm})=(0.5,0.5)$, and we define the set $\Omega_\text{c}$ (\textit{i.e. }the ‘centre of the tumour’) via definition~\eqref{eq:centreTum} as
$$\Omega_\text{c}:=\left\{\mathbf{x_i}\equiv(x_i,y_j)\in [0,1]^2 : \sqrt{(x_i-x_\text{cm})^2+(y_j-y_\text{cm})^2}<0.144 \right\}.$$
The set $\Omega_\text{c}$ corresponds to approximately $65\%$ of the region initially occupied by the tumour and, therefore, the set $\Omega_\text{m} = \Omega \setminus \Omega_\text{c}$ (\textit{i.e. }the ‘margin of the tumour’) comprises the remaining $35\%$ of the tumour region.
%$$\Omega_m:=\left\{x_i\equiv(x_i,y_j)\in \Omega^2, \; \sqrt{(x_i-x_{cm1}^*)^2+(x_j-x_{cm2}^*)^2}<R_2 \right\}\setminus \Omega_c. $$
\paragraph{Parameter values}
Unless otherwise specified, we use the parameter values listed in Table~\ref{table2}. Here, the value of the
parameter $\alpha_n$ is consistent with previous measurement and estimation studies on the
dynamics of tumour cells by~\cite{christophe2015biased}, who calculated the estimated proliferation rate of a tumour cell by using the average duplication time of melanoma cells. The values of the diffusion coefficient $\beta_\phi$ and decay rate $\kappa_\phi$ of the chemoattractant correspond to those used in~\citep{cooper, matzavinos2004mathematical}. Moreover, the range of values of the secretion rate $\alpha_\phi$
is consistent with the values used in~\citep{atsou2020size}. %The values of the parameters $\mu_n$ and $\mu_c$ are chosen so as to ensure that the equilibrium sizes
%of the two cell populations in isolation are biologically relevant. 
To explore a wide range of biological situations corresponding to different degrees of immune infiltration, we use an arbitrary range of values for the parameter $w_{\max}$. Moreover, on the basis of the considerations drawn in~\citep{bubba2020discrete}, we define $\phi_{\max}:= A\,w_{\max}$, where $A \in \mathbb{R^*_+}$ is a scaling factor that ensures unit consistency. We verified via preliminary numerical simulations that, under initial conditions~\eqref{eq:init_chem} and~\eqref{eq:init_chem_cont} for the concentration of the chemoattractant, if $A=1$  then the concentration of chemoattractant remains uniformly smaller than $\phi_{\max}$ for all times. Hence, we set $\phi_{\max} = w_{\max}$ to carry out numerical simulations. Finally, the values of the parameters $\lambda$ and $\nu$ correspond to values of $\beta_c$ and $\gamma_c$ that are consistent with those reported in~\citep{matzavinos2004mathematical} and~\citep{atsou2020size}, respectively. 
Given the values of the parameters $\beta_c$ and $\gamma_c$ chosen to carry out numerical simulations of the continuum model, the following definitions are used for the hybrid model
$$\lambda:=\beta_c\frac{4\tau}{\chi^2}   \quad \text{and} \quad  \nu:=\gamma_c\frac{4\phi_{\max}\tau}{\chi^2}  $$
so that conditions \eqref{eq:conditionsContMod} are met.

\begin{table}[hp!]
	\footnotesize
\caption{Model parameters and related values used in numerical simulations.}
\centering
{\renewcommand{\arraystretch}{1.5}
 	\begin{tabular}{p{1.3cm} p{3.2cm} p{4cm} p{1.5cm}}
  \hline
  Phenotype  & Description & Value \& Units & Reference \\
  \hline  
   \textbf{Domain}  & Space-step in the $x$ and $y$ direction& $\chi= \;0.016$ [\textit{cm}] &  \\  
     &Time-step&  $\tau = 10^{-4}$ [\textit{days}] & \\
      & Final time &   $t_f=15$ [\textit{days}]  & \\ \\
      \textbf{Tumour cells} & Cell density at position $\mathbf{x}$ &$n(\mathbf{x},t)\geq 0$ [\textit{cells/cm$^2$}]&  \\
   & Initial number  &$\rho_n(0) =45228$ [\textit{cells}] &  \\ 
& Proliferation rate  & $\alpha_n=1.5$ [\textit{1/days}] & ~\citep{christophe2015biased}  
 \\
  & Rate of death due to competition between tumour cells  & $\mu_n=1.25\times 10^{-5}$ [\textit{1/(days cells)}] &   
 \\ 
& Level of efficiency of T cells  & $\zeta_n\in[0.001,1]$ [\textit{1/(days cells)}] & \\
& Radius of interaction between tumour cells and T cells  & $\theta=3\times 0.016$ [\textit{cm}] &     
 \\\\
 \textbf{T cells}  & Cell density at position $\mathbf{x}$ & $c(\mathbf{x},t) \geq 0 $ [\textit{cells/cm$^2$}]&  \\
  & Initial number  &$\rho_c(0) =8960$ [\textit{cells}] &  \\ 
  & Prop. const. for influx rate & $\alpha_c=6$ [\textit{cells/(cm$^2$ days mol)}]  &   \\ 
  & Chemotactic sensitivity (hybrid model) &   
   $\nu=\gamma_c\frac{4\phi_{\max}\tau}{\chi^2}$ &   \\ 
     & Chemotactic coefficient (continuum model) & $\gamma_c=10$ [\textit{cm$^2$/(days mol)}]  &~\citep{atsou2020size}   \\ 
     & Random movement prob. (hybrid model) &  
   $\lambda=\beta_c\frac{4\tau}{\chi^2}$ &   \\
     & Diffusion coefficient (continuum model)& $\beta_c=1\times 10^{-3}$ [\textit{cm$^2$/days}] &  ~\citep{matzavinos2004mathematical} \\ 
     & Total cell density above which T cell movement is impaired & $w_{\max}=[0.74\times10^5, 8.88\times10^5]$ [\textit{cells/$cm^2$}]   &   \\ 
   & Rate of death due to competition between T cells  & $ \mu_c=6\times10^{-6}$ [\textit{1/(days cells)}] &  \\   \\
      
     \textbf{Chemoattr.} 
    & Concentration at position $\mathbf{x}$&  $\phi(\mathbf{x},t)\geq 0 $ [\textit{mol/cm}$^2$] &  \\
     & Total amount  &  $\phi_{tot}(t)\geq 0 $ [\textit{mol}] &  \\
     & Diffusion coefficient & $ \beta_\phi=10^{-1}$ [\textit{cm$^2$/day}] &~\citep{matzavinos2004mathematical}  \\
       & Secretion rate & $\alpha_\phi\in[0.001,1.5]$[\textit{mol/(cells days)}]&\citep{atsou2020size}\\
       & Decay rate & $\kappa_\phi=2 $ [\textit{1/days}] &~\citep{cooper}   \\
       & Maximum concentration & $\phi_{\max}=Aw_{\max},$ with $A=1$  [\textit{mol/cm}$^2$]  &   \\
  \hline
\end{tabular}
}
\label{table2}
\end{table}

\subsection{Baseline scenario corresponding to tumour eradication}
\begin{figure}[!t]
\centering
\includegraphics[scale=0.7]{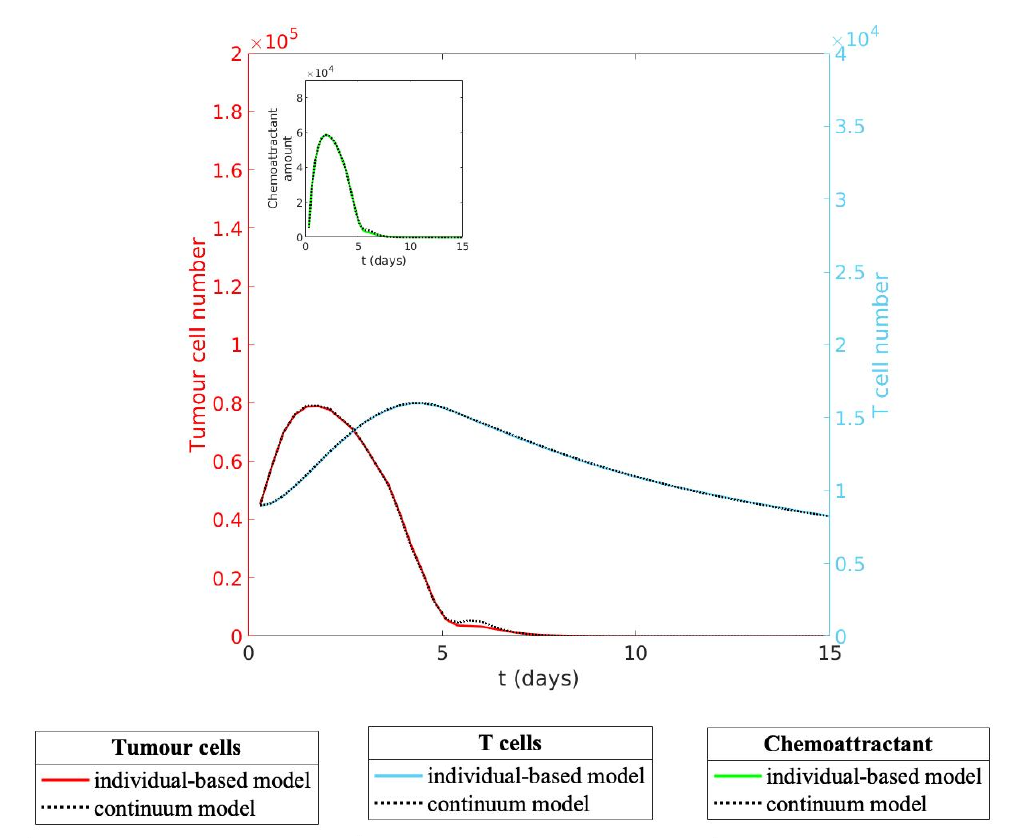}
\caption{\textbf{Baseline scenario: evolution of the numbers of tumour cells and T cells, and the total amount of chemoattractant}. Plots of the time evolution of the number of tumour cells $\rho_n(t)$, the number of T cells $\rho_c(t)$, and the total amount of chemoattractant $\phi_{tot}(t)$ (in the inset) of the hybrid
model (solid, coloured lines) and the continuum model (dotted, black lines) for a choice of parameters that results in the eradication of the tumour. Here, $\zeta_n=0.004$ and all the other parameters are as in Table~\ref{table2} with $\alpha_\phi=1.5$ and $w_{\max}= 2.96\times10^5$. The results from the hybrid model correspond to the average over three simulations and
the related variance is displayed by the coloured areas surrounding the curves.}
\label{fig:evolOverTime}
\end{figure}

\begin{figure}[!t]
\centering
\includegraphics[scale=0.75]{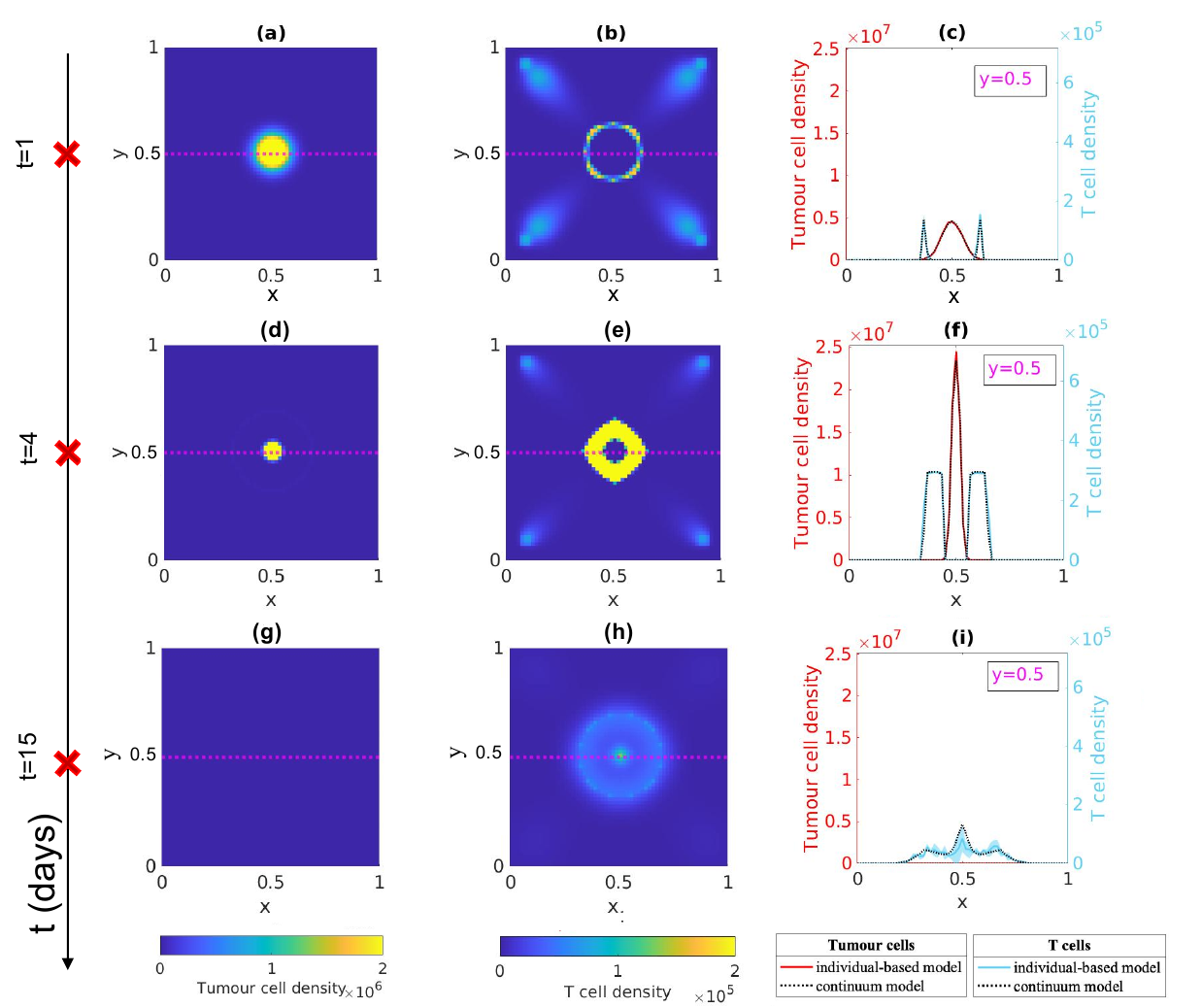}
\caption{\textbf{Baseline scenario: evolution of the spatial distributions of cells over time}. Panels \textbf{(a)-(d)-(g)}
display the plots of the density of tumour cells $n(x,y,t)$ and panels \textbf{(b)-(e)-(h)} display the plots of the density of T cells $c(x,y,t)$ of the continuum model
at three successive time instants – \textit{i.e.} $t=1$, $t=4$, and $t=15$. The pink dashed lines highlight the 1D cross-section corresponding to $y = 0.5$. Panels \textbf{(c)-(f)-(i)} display the corresponding side on view plot of the densities of tumour cells $n(x,y,t)$ and T cells $c(x,y,t)$ of the hybrid model (solid, coloured lines) and continuum model (dotted, black lines) (\textit{i.e.} at $y=0.5$ and $t=1$, $t=4$, and $t=15$). Here, $\zeta_n=0.004$ and all the other parameters are as in Table~\ref{table2} with $\alpha_\phi=1.5$ and $w_{\max}= 2.96\times10^5$. The results from the hybrid model correspond to the average over three simulations and
the related variance is displayed by the coloured areas surrounding the curves.}
\label{fig:evolCellsOverTime}
\end{figure}

As mentioned earlier, we first establish a baseline scenario where the parameter $\zeta_n$ is high enough so that T cells are able to eradicate the tumour. The plots in Fig.~\ref{fig:evolOverTime} and Fig.~\ref{fig:evolCellsOverTime} summarise the results of simulations of the hybrid and continuum models obtained under this scenario.

After initial growth, the number of tumour cells decreases steadily over time until tumour cells are completely eliminated ({\it cf.} Fig.~\ref{fig:evolOverTime}). The chemoattractant produced by tumour cells triggers the inflow of T cells through blood vessels and the movement of T cells towards the tumour ({\it cf.} Fig.~\ref{fig:evolCellsOverTime}). Since the value of $\zeta_n$ is sufficiently large, once T cells are close enough to tumour cells they start eliminating them until eradication ({\it cf.} Fig.~\ref{fig:evolCellsOverTime}). When the number of tumour cells decreases, the total amount of chemoattractant decays as well, thus triggering a reduction in the inflow of T cells, which initiates a decrease in the number of T cells ({\it cf.} Fig.~\ref{fig:evolOverTime}). 

Both Fig.~\ref{fig:evolOverTime} and Fig.~\ref{fig:evolCellsOverTime} indicate that there is an excellent quantitative agreement between numerical solutions of the continuum model \eqref{continuum_model} and the results of numerical simulations of the hybrid model. 
%The only exception is given by the plot in Fig.~\ref{fig:evolCellsOverTime}\textbf{(i)}, which shows that the presence of a lower T cell density may lead to more pronounced demographic stochasticity, and a less regular cell density distribution. In fact, this causes a reduction in the quality of the approximations employed in the formal derivation of the deterministic continuum model from the hybrid model.
\subsection{Emergence of hot, altered and cold tumour scenarios}
\label{Immunoscore and emergence of hot, altered and cold tumour scenarios}
\begin{figure}[!t]
\centering
\includegraphics[scale=0.63]{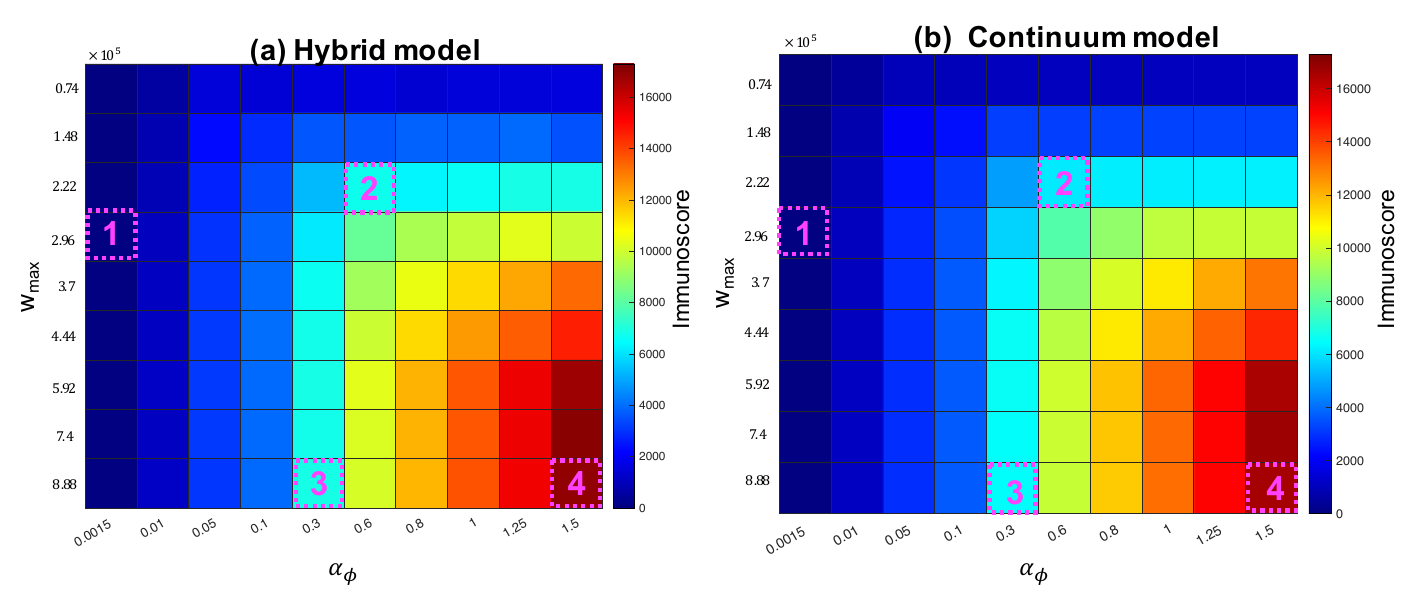}
\caption{\textbf{Immunoscore}. The heat map in panel \textbf{(a)} displays the value of the immunoscore computed via~\eqref{eq:immunoscore_fin} at the end of numerical simulations of the hybrid model for different combinations of $\alpha_\phi$ and $w_{\max}$. For each given value of $\alpha_\phi$ and $w_{\max}$, the values of the other  parameters  are as  in  Table~\ref{table2}, with $\zeta_n=0.00012$. This heat map matches with the corresponding heat map obtained for the continuum model, which is displayed in panel \textbf{(b)}. Sample dynamics of the numbers and densities of tumour cells and T cells for the values of parameters $\alpha_\phi$ and $w_{\max}$ corresponding to the dotted pink squares \textbf{1-4} are displayed in the plots of Fig.~\ref{fig:evolTempsTOT} and Figs.~\ref{fig:evolTempsCold}-\ref{fig:evolTempsHot}, respectively. %Although the specific colors of these areas can vary according to the values of the other parameters of the model, the trends of spatial distribution of T cells and tumour cells in the case of hot, altered-immunosuppressed, altered-excluded and cold tumours remain qualitatively similar to those shown in Figures~\ref{fig:evolTempsHot}-\ref{fig:evolTempsCold}. 
}
\label{fig:immunoscore}
\end{figure}

\begin{figure}[!t]
\centering
\includegraphics[scale=0.7]{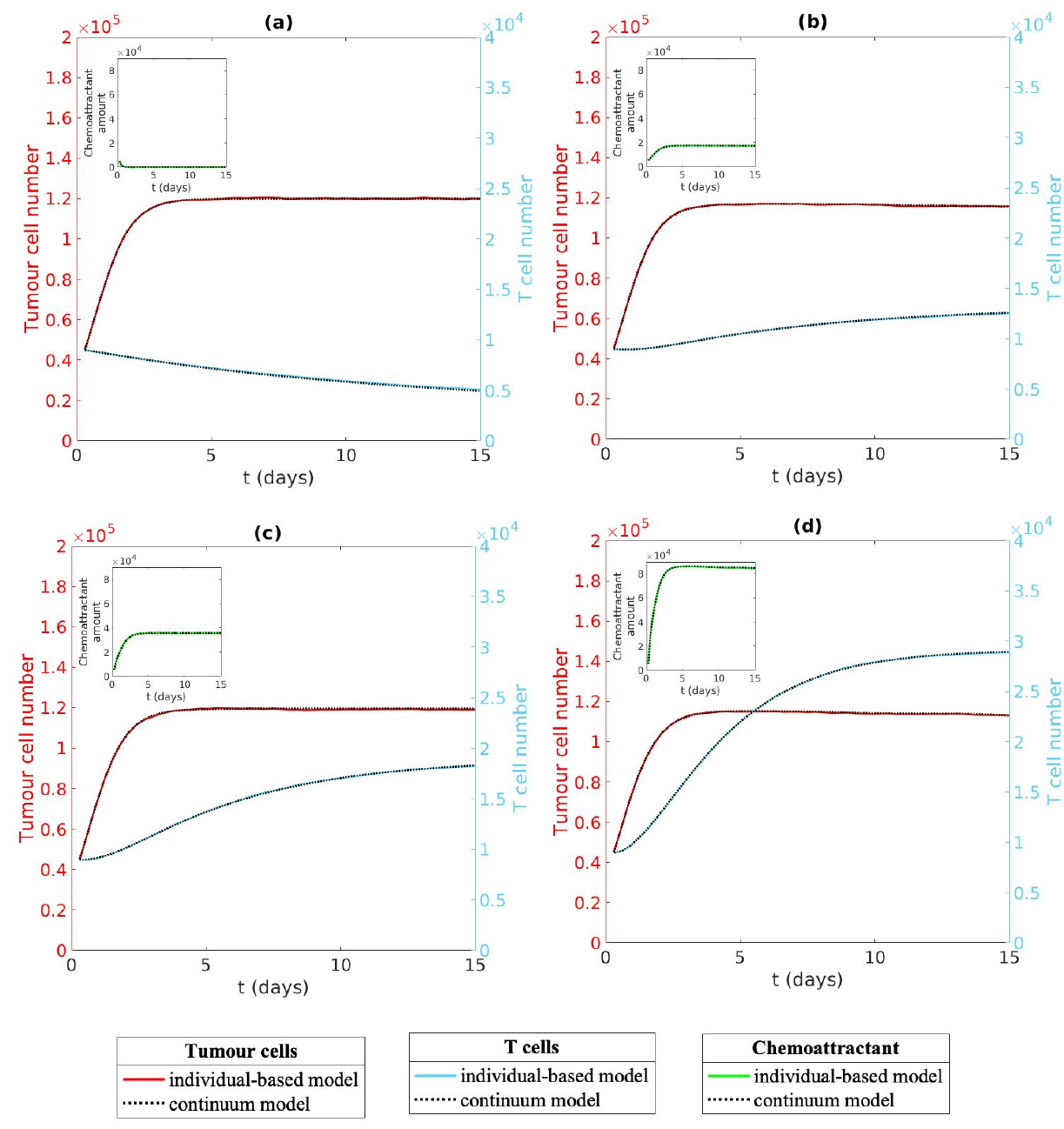}
\caption{\textbf{Sample dynamics of the numbers of tumour cells and T cells in hot, altered-immunosuppressed, altered-excluded and cold tumour scenarios.} Time evolution of the number of tumour cells $\rho_n(t)$, the number of T cells $\rho_c(t)$, and the total amount of chemoattractant $\phi_{tot}(t)$ (in the insets) of the hybrid
model (solid, coloured lines) and the continuum model (dotted, black lines) for values of $\alpha_\phi$ and $w_{\max}$ corresponding to the dotted pink squares \textbf{1-4} of Fig. \ref{fig:immunoscore}. Sufficiently low values of $\alpha_\phi$ lead to the emergence of cold tumour scenarios (panel \textbf{(a)}); intermediate values of $\alpha_\phi$ and sufficiently high values of $w_{\max}$ lead to the emergence of altered-immunosuppressed tumour scenarios (panel \textbf{(b)}); intermediate values of $\alpha_\phi$ and sufficiently small values of $w_{\max}$ lead to the emergence of altered-excluded tumour scenarios (panel \textbf{(c)}); and sufficiently high values of $\alpha_\phi$ and $w_{\max}$ lead to the emergence of hot tumour scenarios (panel \textbf{(d)}). %Four values for the parameter $\alpha_\phi$ are tested: $\alpha_\phi=0.0015$ (panel \textbf{(a)}), $\alpha_\phi=0.1$ (panel \textbf{(b)-(c)}) and $\alpha_\phi=0.0015$ (panel \textbf{(d)}). 
To obtain these results, we used the values of the parameters $\alpha_\phi$ and $w_{\max}$ corresponding to the  dotted pink squares \textbf{1} (panel \textbf{(a)}), \textbf{3} (panel \textbf{(b)}), \textbf{2} (panel \textbf{(c)}), and \textbf{4} (panel \textbf{(d)}) of Fig. \ref{fig:immunoscore}. %: $(\alpha_\phi,w_{\max})=(0.0015, 2.96\times 10^5$) (panel \textbf{(a)}), $(\alpha_\phi,w_{\max})=(0.3,0.75\times 2.96\times 10^5$) (panel \textbf{(b)}), $(\alpha_\phi,w_{\max})=(0.6,3\times 2.96\times 10^5$) (panel \textbf{(c)}) and $(\alpha_\phi,w_{\max})=(1.5,3\times 2.96\times 10^5$) (panel \textbf{(d)}). 
All the other parameters are as in Table~\ref{table2}, with $\zeta_n=0.00012$. %Two values for the parameter $w_{\max}$ are tested: $w_{\max}=6\times 2.96\times 10^5$ (panel \textbf{(a)-(b)}) and $w_{\max}=0.5\times 2.96\times 10^5$ (panel \textbf{(c)-(d)}).  
The results from the hybrid model correspond to the average over three simulations and
the related variance is displayed by the coloured areas surrounding the curves.}
\label{fig:evolTempsTOT}
\end{figure}

\begin{figure}[!t]
\centering
\includegraphics[scale=0.55]{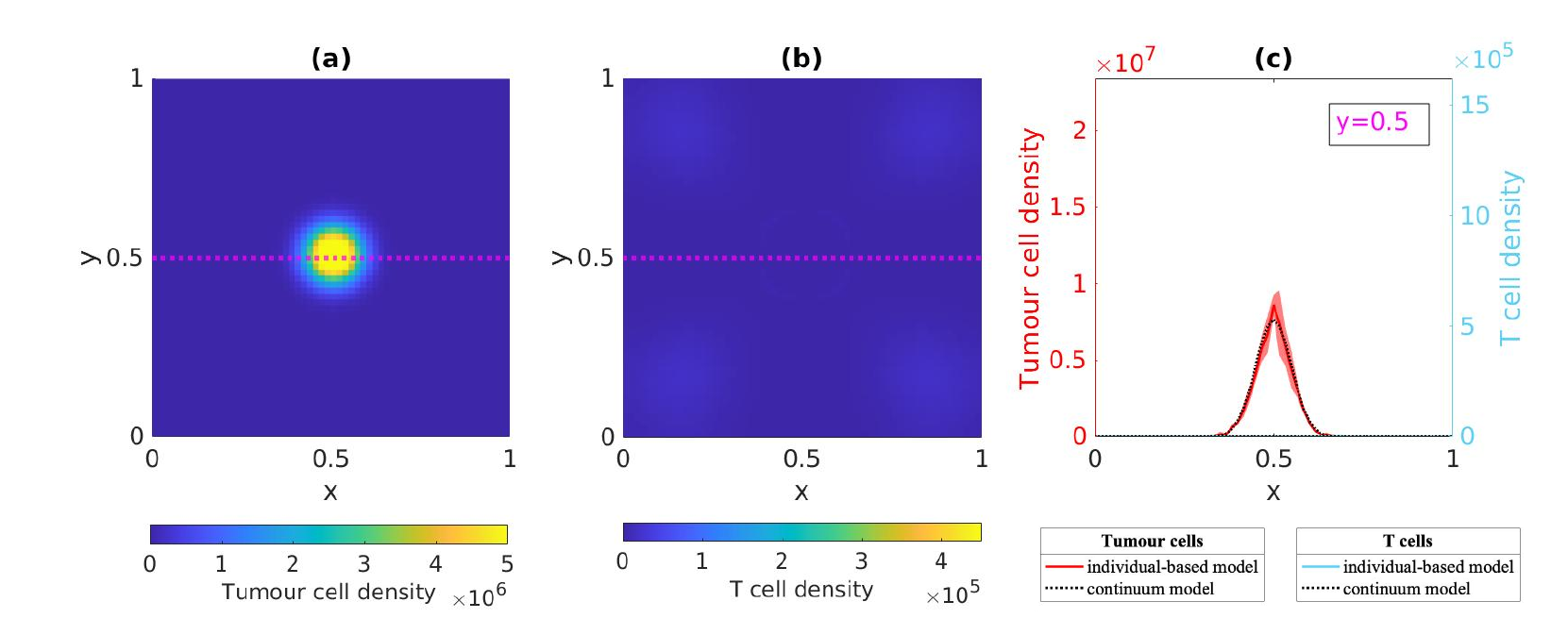}
\caption{\textbf{Sample spatial distributions of cells in cold tumour scenarios}. Panels \textbf{(a)} and \textbf{(b)} display, respectively, the densities of tumour cells $n(x,y,t)$ and T cells $c(x,y,t)$ of the continuum model at the end of numerical simulations (\textit{i.e.} at $t=15$) for a choice of parameters that results in the emergence of a cold tumour scenario ({\it cf.} dotted pink square \textbf{1} in Fig. \ref{fig:immunoscore}). The pink dashed line highlights the 1D cross-section corresponding to $y = 0.5$. Panel \textbf{(c)} displays the corresponding side on view plot of the densities of tumour cells $n(x,y,t)$ and T cells $c(x,y,t)$ of the hybrid model (solid, coloured lines) and continuum model (dotted, black lines) (\textit{i.e.} at $y=0.5$ and $t=15$). Here, $\zeta_n=0.00012$, $\alpha_\phi=0.0015$, $w_{\max}= 2.96\times10^5$, and all the other parameters are as in Table~\ref{table2}. The results from the hybrid model correspond to the average over three simulations and
the related variance is displayed by the coloured areas surrounding the curves.}
\label{fig:evolTempsCold}
\end{figure}

\begin{figure}[!t]
\centering
\includegraphics[scale=0.55]{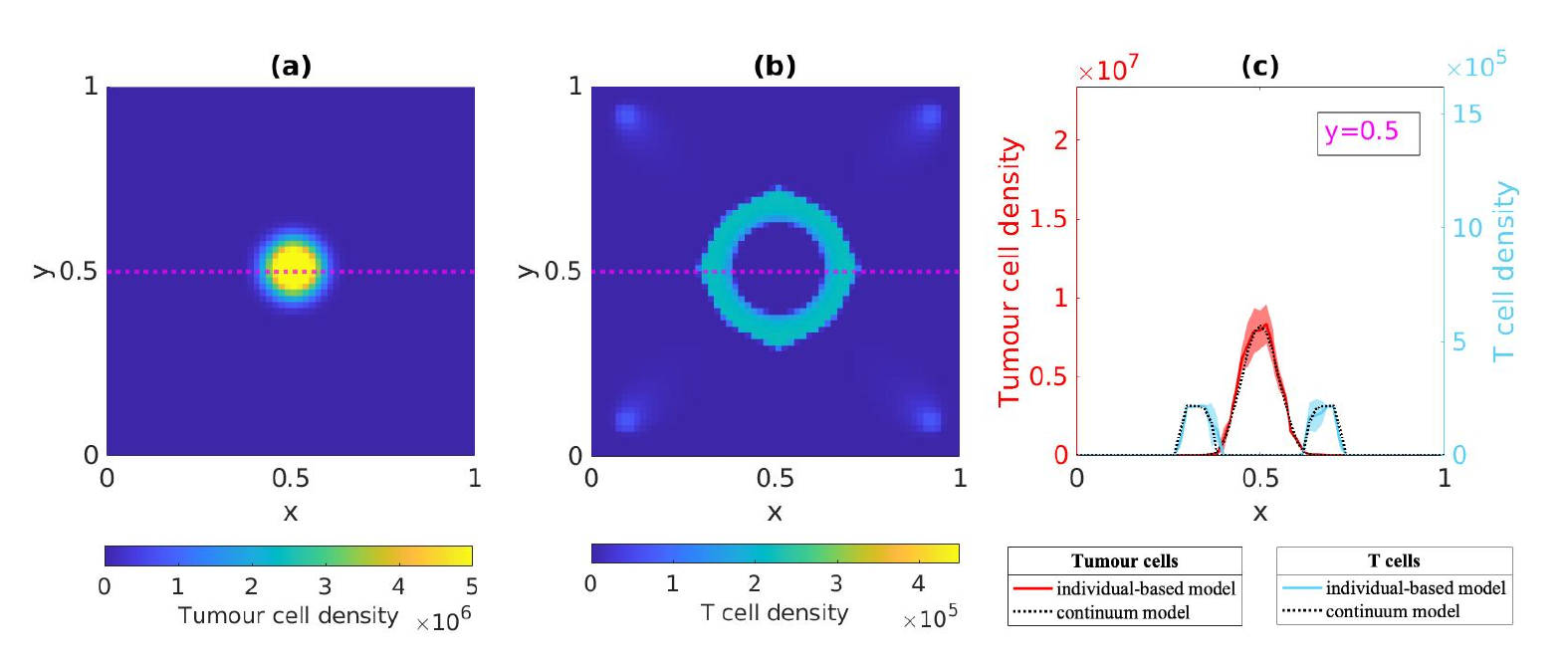}
\caption{\textbf{Sample spatial distributions of cells in altered-excluded tumour scenarios}. Panels \textbf{(a)} and \textbf{(b)} display, respectively, the densities of tumour cells $n(x,y,t)$ and T cells $c(x,y,t)$ of the continuum model at the end of numerical simulations (\textit{i.e.} at $t=15$) for a choice of parameters that results in the emergence of a altered-excluded tumour scenario ({\it cf.} dotted pink square \textbf{2} in Fig. \ref{fig:immunoscore}). The pink dashed line highlights the 1D cross-section corresponding to $y = 0.5$. Panel \textbf{(c)} displays the corresponding side on view plot of the densities of tumour cells $n(x,y,t)$ and T cells $c(x,y,t)$ of the hybrid model (solid, coloured lines) and continuum model (dotted, black lines) (\textit{i.e.} at $y=0.5$ and $t=15$). Here, $\zeta_n=0.00012$, $\alpha_\phi=0.15$, $w_{\max}= 2.22\times10^5$, and all the other parameters are as in Table~\ref{table2}. The results from the hybrid model correspond to the average over three simulations and
the related variance is displayed by the coloured areas surrounding the curves.}
\label{fig:evolTempsExcluded}
\end{figure}

\begin{figure}[!t]
\centering
\includegraphics[scale=0.55]{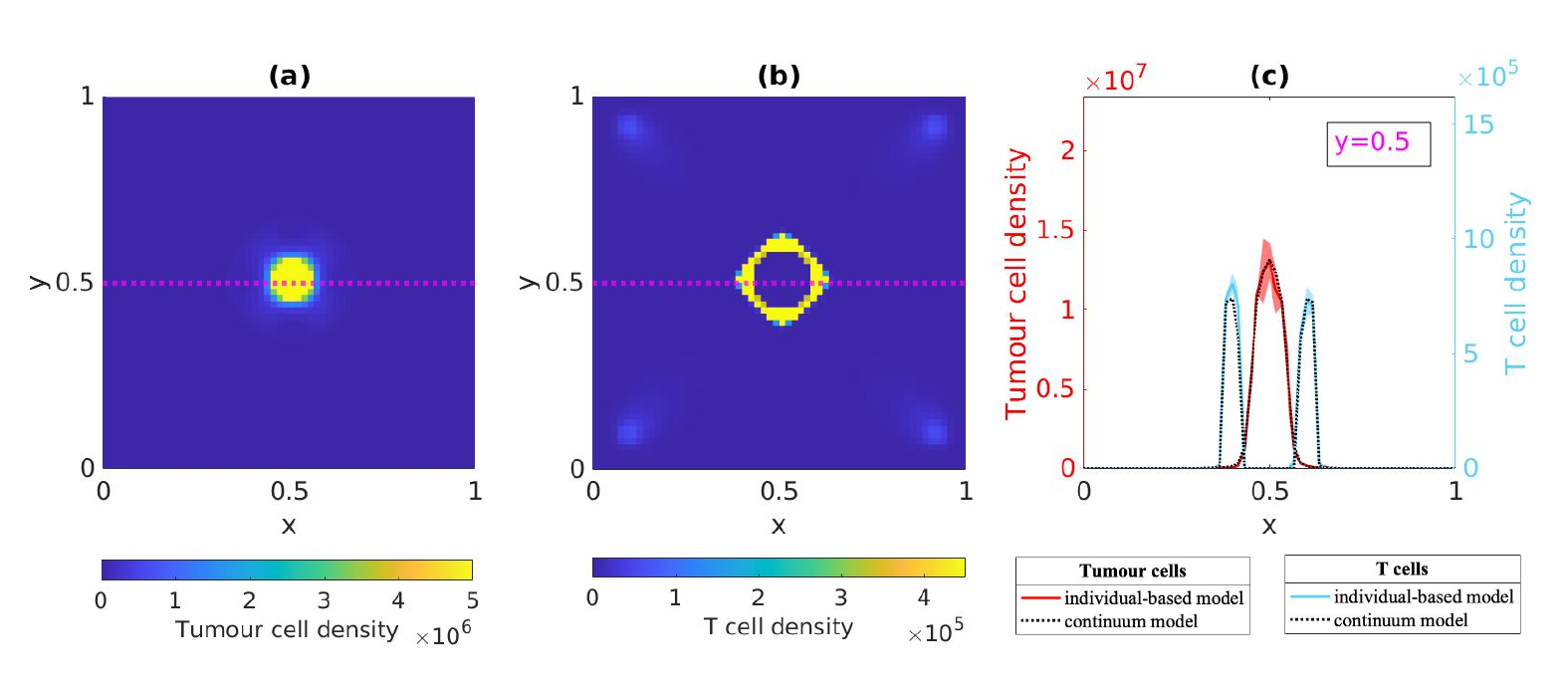}
\caption{\textbf{Sample spatial distributions of cells in altered-immunosuppressed tumour scenarios}. Panels \textbf{(a)} and \textbf{(b)} display, respectively, the densities of tumour cells $n(x,y,t)$ and T cells $c(x,y,t)$ of the continuum model at the end of numerical simulations (\textit{i.e.} at $t=15$) for a choice of parameters that results in the emergence of a altered-immunosuppressed tumour scenario ({\it cf.} dotted pink square \textbf{3} in Fig. \ref{fig:immunoscore}). The pink dashed line highlights the 1D cross-section corresponding to $y = 0.5$. Panel \textbf{(c)} displays the corresponding side on view plot of the densities of tumour cells $n(x,y,t)$ and T cells $c(x,y,t)$ of the hybrid model (solid, coloured lines) and continuum model (dotted, black lines) (\textit{i.e.} at $y=0.5$ and $t=15$). Here, $\zeta_n=0.00012$, $\alpha_\phi=0.15$, $w_{\max}= 8.88\times10^5$, and all the other parameters are as in Table~\ref{table2}. The results from the hybrid model correspond to the average over three simulations and
the related variance is displayed by the coloured areas surrounding the curves.}
\label{fig:evolTempsImmunosuppr}
\end{figure}

\begin{figure}[!t]
\centering
\includegraphics[scale=0.55]{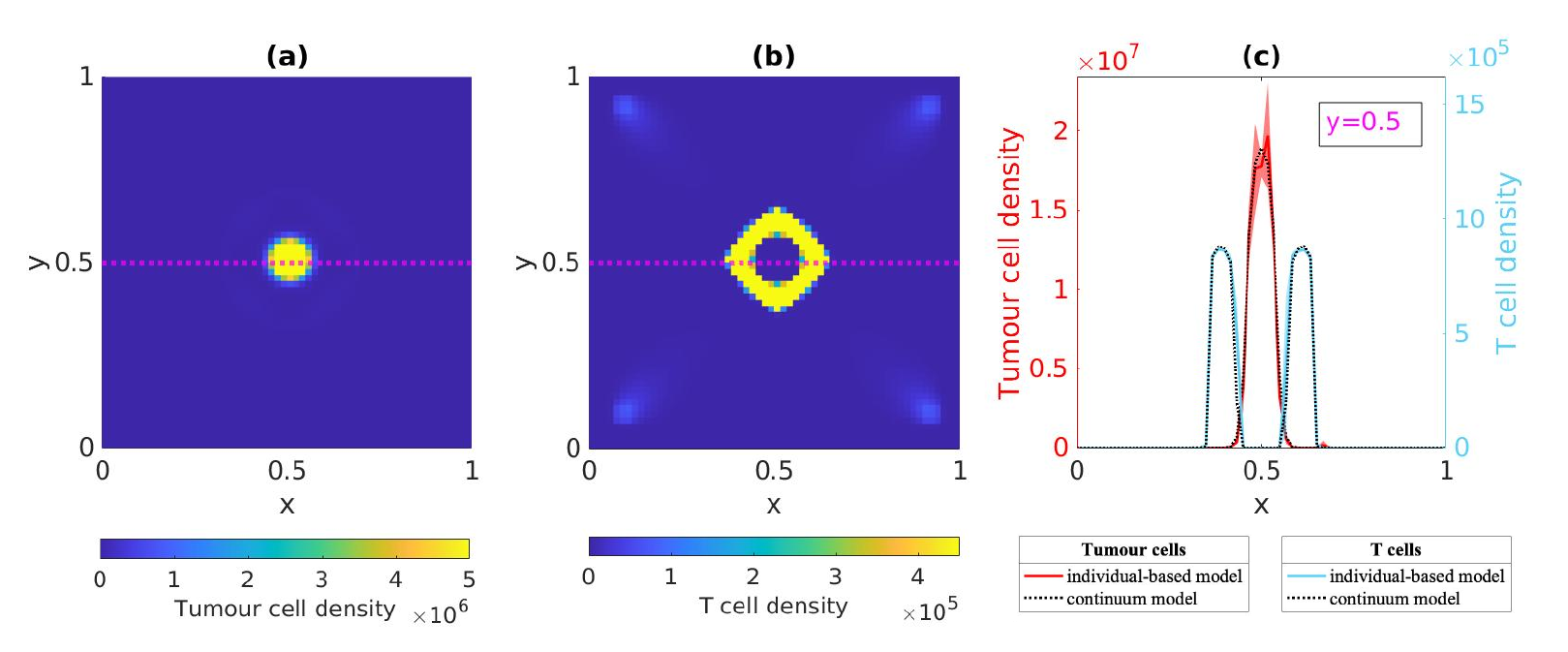}
\caption{\textbf{Sample spatial distributions of cells in hot tumour scenarios}. Panels \textbf{(a)} and \textbf{(b)} display, respectively, the density of tumour cells $n(x,y,t)$ and T cells $c(x,y,t)$ of the continuum model at the end of numerical simulations (\textit{i.e.} at $t=15$) for a choice of parameters that results in the emergence of a hot tumour scenario ({\it cf.} dotted pink square \textbf{4} in Fig. \ref{fig:immunoscore}). The pink dashed line highlights the 1D cross-section corresponding to $y = 0.5$. Panel \textbf{(c)} displays the corresponding side on view plot of the densities of tumour cells $n(x,y,t)$ and T cells $c(x,y,t)$ of the hybrid model (solid, coloured lines) and continuum model (dotted, black lines) (\textit{i.e.} at $y=0.5$ and $t=15$). Here, $\zeta_n=0.00012$, $\alpha_\phi=1.5$, $w_{\max}= 8.88\times10^5$, and all the other parameters are as in Table~\ref{table2}. The results from the hybrid model correspond to the average over three simulations and
the related variance is displayed by the coloured areas surrounding the curves.}
\label{fig:evolTempsHot}
\end{figure}

We now consider a lower value of the parameter $\zeta_n$ in order to explore biological scenarios in which the cytotoxic action of T cells is less effective, for example due to high expression of PD1 inhibitory receptors and PD-L1 ligands on the surface of T cells and tumour cells. As mentioned earlier, we wish to investigate how the spatial distribution of T cells within the tumour varies depending on the value of the parameters $\alpha_\phi$ and $w_{\max}$. Therefore, we perform numerical simulations holding all parameters constant but considering different combinations of $\alpha_\phi$ and $w_{\max}$. For each pair of values considered, we stored the resulting dynamics of the densities of tumour cells and T cells along with the dynamics of the corresponding cell numbers, and the final value of the immunoscore computed via~\eqref{eq:immunoscore_fin}. The results obtained are summarised by the heat maps in Fig.~\ref{fig:immunoscore} and the plots in Figs.~\ref{fig:evolTempsTOT}-\ref{fig:evolTempsHot}.
\paragraph{Low immunoscore and emergence of cold tumour scenarios}
As shown by the blue regions on the left side of the two heat maps of Fig.~\ref{fig:immunoscore}, for sufficiently small values of $\alpha_\phi$, the immunoscore is relatively low independently of the value of $w_{\max}$. This is due to the small concentration of chemoattractant present in the domain, which poses limitations to the inflow and movement of T cells towards the tumour. %In this case, the value of the immunoscore defined
%via \eqref{eq:immunoscore_fin} at the final time of our simulations is relatively small. 
Hence, in the framework of our model, this parameter range corresponds to the emergence of cold tumour scenarios. 

Sample dynamics of the numbers and densities of tumour cells and T cells for the values of $\alpha_\phi$ and $w_{\max}$ corresponding to the dotted pink square \textbf{1} in Fig.~\ref{fig:immunoscore} are displayed in the plots in Fig.~\ref{fig:evolTempsTOT}\textbf{(a)} and Fig.~\ref{fig:evolTempsCold}. As shown by Fig.~\ref{fig:evolTempsTOT}\textbf{(a)}, for sufficiently small values of $\alpha_\phi$, the total amount of chemoattractant in the domain is too small to trigger a sufficiently high inflow of T cells that can compensate for the loss caused by T cell death. As a result, the number of T cells decreases over time. Moreover, there is a shallow gradient of the chemoattractant, which results in a slow movement of T cells towards the tumour. As a result, as shown by Fig.~\ref{fig:evolTempsCold}\textbf{(b)-(c)}, at the end of simulations, the density of T cells around the tumour is almost zero and T cells are still very much concentrated in the proximity of the blood vessels ({\it i.e.} their entry points).
\paragraph{Intermediate immunoscore and emergence of altered tumour scenarios} The light blue regions of the heat maps of Fig. \ref{fig:immunoscore} indicate that there are two possible parameter ranges giving rise to an intermediate immunoscore. The first one corresponds to intermediate values of $\alpha_\phi$ along with intermediate to large values of $w_{\max}$, while the second one corresponds to larger values of $\alpha_\phi$ along with small values of $w_{\max}$. In the framework of our model, altered tumour scenarios emerge under these parameter ranges.%, which can then be further classified as altered-immunosuppressed or
%altered-excluded based on the distribution of T cells at the centre and margin of the tumour. 

Sample dynamics of the numbers and densities of tumour cells and T cells for the values of $\alpha_\phi$ and $w_{\max}$ corresponding to the dotted pink squares \textbf{2} and \textbf{3} in Fig.~\ref{fig:immunoscore} are displayed in the plots in Fig.~\ref{fig:evolTempsTOT}\textbf{(b)}-\textbf{(c)}, Fig.~\ref{fig:evolTempsExcluded} and Fig.~\ref{fig:evolTempsImmunosuppr}. The results of Fig.~\ref{fig:evolTempsTOT}\textbf{(b)}-\textbf{(c)} show that increasing the value of $\alpha_\phi$ leads to a progressive increase in the total amount of chemoattractant.
This in turn results in an increased inflow of T cells and facilitates the movement of T cells towards the tumour. The spatial distribution of T cells within the tumour varies depending on the value of $w_{\max}$. Fig.~\ref{fig:evolTempsExcluded} shows that smaller values of $w_{\max}$  lead  to  an  accumulation  of  T  cells  at the margin of the tumour, which corresponds to an altered-excluded tumour scenario. On the other hand, larger values of $w_{\max}$ promote the infiltration of T cells into the tumour and lead to an altered-immunosuppressed tumour scenario (see Fig.~\ref{fig:evolTempsImmunosuppr}). %These different outcomes are due to the different choice of values of $w_{\max}$, which regulate the maximum attainable infiltration level of T cells within the tumour.  
\paragraph{High immunoscore and emergence of hot tumour scenarios} Finally, as shown by the red regions on the bottom-right side of Fig.~\ref{fig:immunoscore}, for large values of $\alpha_\phi$ and $w_{\max}$, the value of the immunoscore is relatively high. In the framework of our model, this parameter range corresponds to the emergence of hot tumour scenarios.  

Sample dynamics of the numbers and densities of tumour cells and T cells for the values of $\alpha_\phi$ and $w_{\max}$ corresponding to the dotted pink square \textbf{4} in Fig.~\ref{fig:immunoscore} are displayed in the plots in Fig.~\ref{fig:evolTempsTOT}\textbf{(d)} and Fig.~\ref{fig:evolTempsHot}. When $\alpha_\phi$ is high enough, the larger amount of chemoattractant promotes the inflow of a larger number of T cells (see Fig.~\ref{fig:evolTempsTOT}\textbf{(d)}). Moreover, Fig.~\ref{fig:evolTempsHot} shows that, similarly to the altered-immunosuppressed tumour scenario, larger values of $w_{\max}$ facilitate the infiltration of T cells into the tumour. As the number of infiltrated T cells is larger than in the previous scenarios, the immune action is slightly more efficient and thus leads to a slightly decreased number of tumour cells (see Fig.~\ref{fig:evolTempsTOT}).

\begin{remark}
Although the specific colours of the regions of the heat maps in Fig. \ref{fig:immunoscore} can vary according to the values of the other parameters of the model, the behaviours of the spatial distributions of T cells and tumour cells in the case of hot, altered-immunosuppressed, altered-excluded and cold tumour scenarios remain qualitatively similar to those shown in Figs.~\ref{fig:evolTempsCold}-\ref{fig:evolTempsHot}. Moreover, the heat maps in Fig. \ref{fig:immunoscore}, as well as the plots in Fig.~\ref{fig:evolTempsTOT} and Figs.~\ref{fig:evolTempsCold}\textbf{(c)}-\ref{fig:evolTempsHot}\textbf{(c)} demonstrate that there is an excellent agreement between numerical simulations of the hybrid
and continuum models. This testifies to the robustness of the computational results presented here and the biological insight that they provide.
\end{remark}

\subsection{Immunotheraphy effects}
\label{Immunotheraphy effects}
%tumour cells are not eradicated (\textit{i.e.} when the value of $\zeta_n$ is sufficiently small) and different combinations of the
The results presented in the previous subsection summarise how scenarios corresponding to different levels of T-cell infiltration into the tumour can emerge under different combinations of the values of the parameters $\alpha_\phi$ and $w_{\max}$. We now investigate possible outcomes of immunotheraphy in these different scenarios.

\begin{figure}[!t]
\centering
\includegraphics[scale=0.72]{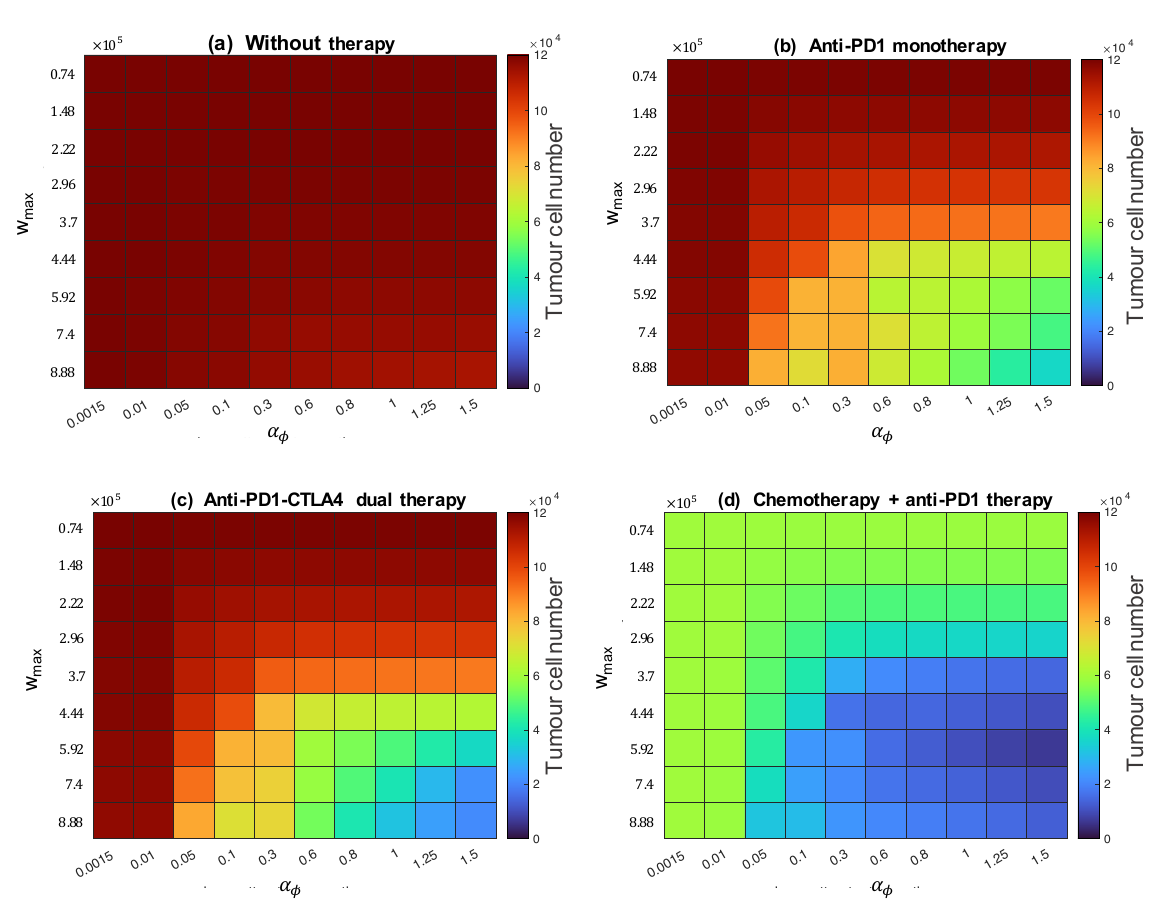}
\caption{\textbf{Immunotherapy effects}. Comparison between the numbers of tumour cells $\rho_n(t)$ at the end of simulations (\textit{i.e.} at $t=10$) of the continuum model for different values of $\alpha_\phi$ and $w_{\max}$ without therapy (panel \textbf{(a)}) and when the effects of different immunotherapies are incorporated by considering different values of the parameters $\zeta_n$, $\alpha_c$ and $\alpha_n$ (panel \textbf{(b)-(d)}. In panel \textbf{(a)}, $\zeta_n=0.00012$ and the values of the other  parameters  are as  in  Table~\ref{table2}. In panel \textbf{(b)}, $\zeta_n=0.001$  and the values of the other  parameters  are as  in  Table~\ref{table2}. In panel \textbf{(c)}, $\zeta_n=0.001$, $\alpha_c=12$  and the values of the other  parameters  are as  in  Table~\ref{table2}. In panel \textbf{(d)}, $\zeta_n=0.001$, $\alpha_c=12$, $\alpha_n=0.75$  and the values of the other  parameters  are as  in  Table~\ref{table2}.}
\label{fig:immunotherapy}
\end{figure}

In order to do this, we consider the same parameter settings used for the numerical simulations of Fig. \ref{fig:immunoscore}, but we allow the level of efficiency of T cells at eliminating tumour cells to be higher (\textit{i.e.} we increase the value of $\zeta_n$). This corresponds to a biological scenario in which the tumour is treated with anti-PD1 monotherapy, which restores immune efficacy~\citep{topalian2012safety}. We also investigate the effects of coupling anti-PD1 therapy with two other therapies. %This is done by changing the value of other parameters of our model that can be associated to therapy effects.
First we explore the effects of anti-PD1 therapy in combination with another immune checkpoint therapy, {\it i.e.} the anti CTLA-4 therapy~\citep{van2015genomic}. To do so, we perform numerical simulations defining
all parameters as in the case of the anti-PD1 therapy but increasing the influx rate of T  cells through blood vessels (\textit{i.e.} the value of the parameter $\alpha_c$). Then, we explore the effects of combining anti-PD1 therapy with chemotherapy, which inhibits tumour cell division, inflames the TME with tumour antigens, and boosts the activation of T cells~\citep{galon2019approaches}. To do so, we perform numerical simulations defining
all parameters as in the scenario of the anti-PD1 therapy but decreasing the proliferation rate of tumour cells (\textit{i.e.} the value of the parameter $\alpha_n$) and increasing the influx rate of T cells through
blood vessels (\textit{i.e.} the value of the parameter $\alpha_c$). The results obtained are displayed in Fig. \ref{fig:immunotherapy}, which shows a comparison between the numbers of tumour cells at the end of numerical simulation in the scenario ‘‘without treatment’’ (\textit{i.e.} with the parameter values considered in Section \ref{Immunoscore and emergence of hot, altered and cold tumour scenarios}) and the three aforementioned scenarios in which the effects of different therapeutic protocols are considered. 

Exploiting the excellent quantitative agreement between the results of numerical simulations of the hybrid and continuum models presented in the previous subsections, here we carry out the numerical simulations of the continuum model only, since they require computational times much smaller than those that would be required by the numerical exploration of the corresponding hybrid model. To obtain the results presented in this subsection, we carried out numerical simulations by using a final time corresponding to 10 days (\emph{i.e.} $t_f=10$). 
\paragraph{Anti-PD1 monotherapy}
Fig. \ref{fig:immunotherapy}\textbf{(b)} displays the number of tumour cells at the end of numerical simulations of the continuum model for parameter settings corresponding to anti-PD1 monotherapy (\textit{i.e.} when only the value of $\zeta_n$ in increased). Comparing these results with those displayed in Fig. \ref{fig:immunotherapy}\textbf{(a)}, we see that, in general, for the same values of parameters $\alpha_\phi$ and $w_{\max}$, increasing the value of  $\zeta_n$ leads to a decrease in the number of tumour cells at the end of simulations. However, when the value of $\alpha_\phi$ is too small (\textit{i.e.} in cold tumour scenarios) or when the value of $w_{\max}$ is too small (\textit{i.e.} in altered-excluded tumour scenarios), increasing $\zeta_n$ has no benefit on the action of T cells against tumour cells. Finally, when the values of $\alpha_\phi$ and $w_{\max}$ are sufficiently large (\textit{i.e.} in hot tumour scenarios) anti-PD1 monotherapy is more effective.  
\paragraph{Anti-PD1-CTLA4 dual therapy}
Fig. \ref{fig:immunotherapy}\textbf{(c)} displays the number of tumour cells at the end of numerical simulations of the continuum model for parameter settings corresponding to anti-PD1-CTLA4 dual therapy (\textit{i.e.} when both the value of $\zeta_n$ and the value of $\alpha_c$ are increased). Comparing these results with those displayed in Fig. \ref{fig:immunotherapy}\textbf{(b)}, we see that increasing the value of $\alpha_c$ along with the value of $\zeta_n$ improves immune efficacy only when the values of $\alpha_\phi$ and $w_{\max}$ are large enough (\textit{i.e.} in hot tumour scenarios). Moreover, for intermediate values of $\alpha_\phi$ (\textit{i.e.} in altered-immunosuppressed tumour scenarios), increasing the value of $\alpha_c$ slightly decreases the number of tumour cells at the end of simulations. Finally, when the values of $\alpha_\phi$ or $w_{\max}$ are too small, increasing $\alpha_c$ has no benefit on the action of T cells against tumour cells. 
\paragraph{Chemotherapy combined with anti-PD1 therapy}
Fig. \ref{fig:immunotherapy}\textbf{(d)} displays the number of tumour cells at the end of numerical simulations of the continuum model for parameter settings corresponding to chemotherapy in combination with anti-PD1 therapy (\textit{i.e.} when the value of $\alpha_n$ is decreased and the values of $\zeta_n$ and $\alpha_c$ are increased). Compared to the other heat maps of Fig. \ref{fig:immunotherapy}, these results show that the number of tumour cells decreases even for small values of $\alpha_\phi$ or $w_{\max}$. Moreover, when the values of $\alpha_\phi$ or $w_{\max}$ are small (\textit{i.e.} in cold and altered-excluded tumour scenarios), the numbers of tumour cells at the end of simulations are similar. % and the lower number of tumour cells is due to the decrease of their proliferation rate $\alpha_n$.
On the other hand, from intermediate to large values of $\alpha_\phi$ and $w_{\max}$, the numbers of tumour cells at the end of simulation decrease as the values of these two parameters increase. % and immune efficiency improves due to the combination of the decrease of $\alpha_n$ and the increase of $\alpha_c$. 
As expected, the larger the values of $\alpha_\phi$ and $w_{\max}$ (\textit{i.e.} the ‘‘hotter" the tumour scenario considered), the more effective the combined action of chemotherapy and anti-PD1 therapy.

\newpage
\section{Discussion and research perspectives}
\label{sec:Conclusive discussion and research perspectives}
\subsection{Discussion}
%In this paper, we have developed a hybrid discrete-continuum modelling framework for the interactions between tumour cells and cytotoxic T cells, which explicitly takes into account the spatial dynamics leading to infiltration of T cells into the tumour. We used a discrete mass-balance equation for the evolution of a chemoattractant, which drives the movement of T cells towards the tumour, while cell dynamics were described by stochastic individual-based models. From our hybrid model, a continuum model comprising a PDE for the T cell density, an IDE for the tumour cell density and a PDE for the concentration of chemoattractant has been formally derived. %Through comparison of both models we found that the results of computational simulations of the hybrid model faithfully mirror the qualitative properties of the solutions to the corresponding continuum model. The results of computational simulations of the hybrid model, which are in excellent agreement with numerical solutions of the continuum model, shed light on the way in which different parameters affect the spatial distribution of T cells within the tumour and allow us to assess the impact of T cell infiltration on the immune response against tumours.

The results that we have presented demonstrate that the level of efficiency of T cells at eliminating tumour cells (\textit{i.e.} the parameter $\zeta_{n}$) plays a key role in tumour-immune competition. In fact, when the value of $\zeta_n$ is large enough, our results indicate that tumour eradication can occur, while lower values of $\zeta_n$ may result in tumour cell survival. This is consistent with experimental and clinical data which point to a key role of immune check-points in immunosuppressing T-cell responses. In fact, the presence of immunosuppressive components in the TME, such as PD1 inhibitory receptors and PD-L1  ligands, decreases the efficiency of T cells at eliminating tumour cells, and can ultimately result in tumour escape~\citep{huang2017t, tumeh2014pd}.%T cells express PD1 after being activated as a mechanism of retro-control; using anti-PD1 antibodies restores the activation of the cytotixicity of these cells and a greater clearance of tumours occurs \. 

Moreover, our numerical results indicate that
when tumour eradication does not occur (\textit{i.e.} when the value of $\zeta_{n}$ is sufficiently small), the secretion rate of the chemoattractant by tumour cells (\textit{i.e.} the parameter $\alpha_\phi$) and the threshold value of the total cell density above which T cell movement is impaired (\textit{i.e.} the parameter $w_{\max}$) have a strong impact on the level of infiltration of T cells into the tumour, and different combinations of the values of these parameters bring about the emergence of four immune-based tumour scenarios. Hot tumour scenarios emerge for high values of $\alpha_\phi$ and $w_{\max}$, and are characterised by a large number of T cells in the centre of the tumour. By displaying a high degree of T cell infiltration, these tumour scenarios provide a fertile ground for immune checkpoint therapies. 
%This is also demonstrated by the numerical results presented in Section \ref{Immunotheraphy effects}.
Altered tumour scenarios emerge for intermediate values of $\alpha_\phi$, and reflect the intrinsic ability of the immune system to effectively mount a T-cell-mediated immune response and the ability of tumour cells to partially
escape such a response. This can either be due to an insufficient number of infiltrated T cells (the immunosuppressed tumour scenarios, which emerge for intermediate to large values of $w_{\max}$) or to the presence of physical barriers that hinder T cell infiltration (the excluded tumour scenarios, which emerge for small values of $w_{\max}$). Finally, cold tumour scenarios emerge for sufficiently small values of $\alpha_\phi$ and $w_{\max}$. These tumour scenarios are characterised by an insufficient
number of T cells both in the centre of the tumour and at its margin, and are invariably associated with poor prognosis. 

We also explored how the outcomes of different immunotherapy protocols can vary in such immune-based tumour scenarios. In particular, our results suggest that increasing the level of efficiency of T cells (\textit{i.e.} the value of the parameter $\zeta_{n}$), which is associated to the effects of anti-PD1 monotherapy, is not sufficient for treating all types of tumour scenarios, and it is particularly ineffective in altered-excluded and cold tumour scenarios. This finding is coherent with experimental observations indicating that anti-PD1 monotherapy is effective only in the context of hot or altered-immunosuppressed tumours, as a certain number of T cells is already infiltrated into the tumour~\citep{galon2019approaches}.

Moreover, the results of our model indicate that, in these two categories of tumours, increasing both  the level of efficiency of T cells and their influx rate (\textit{i.e.} the value of the parameters $\zeta_{n}$ and $\alpha_c$), which are associated with the combined effects of anti-PD1 and anti-CTLA4 therapy, may lead to a better therapeutic outcome. This conclusion is also supported by experimental work showing that anti-PD1-CTLA4 dual therapy may be successful in treating advanced-stage melanoma~\citep{wolchok2017overall}, renal-cell carcinoma~\citep{motzer2018nivolumab} and non-small-cell lung cancer (NSCLC)~\citep{hellmann2018nivolumab}, resulting in regulatory approval. However, our results suggest that prognosis in altered-excluded and cold tumour scenarios may not benefit from the combined effects of these two immune checkpoint inhibitors. Nevertheless, our results indicate that therapeutic strategies promoting the infiltration of T cells could turn altered-excluded tumours into altered-immunosuppressed or hot tumours, helping to decrease the resistance of tumours to the combination of anti-PD1 and anti-CTLA4 therapy. This finding is coherent with experimental observations suggesting that a synergistic effect can be achieved by combining anti-angiogenic therapies, which act on vascular abnormalities facilitating T-cell infiltration,  with immune checkpoint therapies~\citep{tian2017mutual}.

Finally, the outputs of our model suggest that increasing both the level of efficiency of T cells and their influx rate through blood vessels (\textit{i.e.} the values of the parameters $\zeta_{n}$ and $\alpha_c$) and decreasing the proliferation rate of tumour cells (\textit{i.e.} the parameter $\alpha_n$), which may represent the combination of anti-PD1 therapy with chemotherapy, a stronger immune response may be induced. In fact, a proposed approach to overcome the lack of a pre-existing immune response consists in combining a priming therapy that enhances T cell responses (such as chemotherapy) with the removal of co-inhibitory signals (through approaches such as immune checkpoint therapies)~\citep{galon2019approaches}. For example, the success of the combination of anti-PD1 therapy with chemotherapy in metastatic NSCLC has demonstrated the strength of this dual approach~\citep{gandhi2018pembrolizumab}.

%Taken together, the outcomes of our model underline the pivotal role of T cells against tumours and the importance of careful assessment of the pre-existing T cell landscape in the tumour microenvironment. This would make it possible to differentiate specific cases at the moment of diagnosis and to build a solid classification strategy supporting subsequent therapy. Moreover, our results recapitulate the conclusions of previous experimental work showing the importance of personalized tumour treatment and the involvement of combinations with immunotherapies to achieve maximal efficiency~\citep{angell2013immune, galon2006type,kato2017prospects}.\\
%The excellent agreement between the results of numerical simulations of the individual-based and
%continuum models testifies to the robustness of the biological insight gained in this work.
%We also showed that possible differences between cell dynamics produced by the individual-based and
%continuum models can emerge under parameter settings that correspond to sufficiently low cell numbers and more pronounced demographic stochasticity. In fact, these cause a reduction in the
%quality of the approximations employed in the formal derivation of the deterministic continuum model
%from the hybrid model. This demonstrates the importance of integrating individual-based and
%continuum approaches when considering mathematical models for tumour-immune interactions.
\subsection{Research perspectives} We conclude with an outlook on possible extensions of the present work. While here we focused on the role of the secretion rate of the chemoattractant by tumour cells and the threshold value of the total cell density above which T cell movement is impaired, it would be interesting to investigate how other model parameters (\textit{e.g.} the chemotactic sensitivity of T cells) may affect the level of infiltration of T cells into the tumour. Carrying out a more extensive exploration of the model parameter space would ultimately allow more robust biological conclusions to be drawn.
%\hl{This paper is intended to present a proof of concept for the ideas underlying the modelling framework, rather than present a systematic investigation of the conditions on the model parameters leading to the emergence of the four immune-based categories of tumours.
%However, it would be interesting to determine a broader list of parameters of interest and use global sensitivity analysis to study how others biological mechanisms may affect the results of numerical simulations. Moreover, in order to disentangle and quantify the impact of different parameters on the emergence of different levels of infiltration
%of T cells into the tumour, it
%would be useful to integrate numerical simulations of the hybrid model with analytical results of the continuum model. This would enable a more extensive exploration of the model parameter space, which would ultimately allow more robust conclusions to be drawn. }

Moreover, our hybrid modelling framework for the spatial dynamics  of tumour cells and cytotoxic T cells, along with the formal derivation of the corresponding continuum  model,  can be developed further in several ways. For instance, a key factor of the immune response is that T cells express a
unique repertoire of T cell receptors (TCRs)~\citep{coulie2014tumour}, and are capable of detecting and eliminating
tumour cells by recognising specific cancer-associated antigens. The model presented here does not include this aspect, but it could easily be 
extended to do so by introducing, for instance, a variable representing the antigens expressed by tumour cells and the TCR expressed by T cells. This would make it possible to take explicitly into account the effects of both spatial and antigen-specific interactions between tumour cells and T cells, as similarly done in~\citep{kather2017silico, leschiera2022mathematical, macfarlane2018modelling, macfarlane2019stochastic}, and then study the effects of antigen presentation or intra-tumour heterogeneity on immune surveillance.

Only a simplified representation of the action of different types of immunotherapy was considered in this work, but it would be important to carry out a more detailed study of the impact of T-cell infiltration on the dynamics of tumour cells under different immunotherapeutic protocols. In particular, by using optimal control methods for the continuum model, we could investigate the best delivery schedule  of therapeutic agents (\textit{i.e.} the best delivery times and dosages) that make it possible to minimise the number of tumour cells at the end of the treatment and achieve the best therapeutic outcomes~\citep{jarrett2020optimal}. These are all lines of research that we will be pursuing in the near future. 

\section*{Declarations}
\subsection*{Data Availability} 
The datasets generated and analysed during the current study are available from the corresponding author on reasonable request.
\subsection*{Funding} 
E.L. has received funding from the European Research Council (ERC) under the European Union's Horizon2020 research and innovation programme (grant agreement No 740623). \\
T.L. gratefully acknowledges support from the Italian Ministry of University and Research (MUR) through the grant ``Dipartimenti di Eccellenza 2018-2022'' (Project no. E11G18000350001) and the PRIN 2020 project (No. 2020JLWP23) ``Integrated Mathematical Approaches to Socio–Epidemiological Dynamics'' (CUP: E15F21005420006). \\
L.A., E.L. and T.L. gratefully acknowledge support from the CNRS International Research Project ``Modélisation de la biomécanique cellulaire et tissulaire'' (MOCETIBI).
\subsection*{Conflicts of interest} 
The authors declare that they have no conflict of interest.

	\begin{appendices}
\section{Formal derivation of the continuum model}
\label{Formal derivation of the continuum model corresponding to the hybrid model}
Building on the methods employed in~\citep{bubba2020discrete}, we carry out a formal derivation of the deterministic continuum model given by the IDE-PDE-PDE system \eqref{continuum_model} for $d = 1$. Similar methods can be used in the case where $d = 2$.
\subsection{Formal derivation of the IDE for the density of tumour cells $n(x,t)$}	

When tumour cell dynamics are governed by the rules described in Sections~\ref{Tumour cell proliferation and natural death} and~\ref{Effect of the immune system on the tumour}, considering $i \in [0,\mathcal{N}]$, 
between time-steps $k$ and $k+1$ the principle of mass balance gives the following difference equation for the tumour cell density $n_{i}^{k}$:
\begin{equation}
n_{i}^{k+1}=\left[2 \, \tau \alpha_n + 1-\tau(\alpha_n+\zeta_n K_{i}^{k}+\mu_n\rho_n^k \right]n_{i}^{k}.
\label{princ_mass}
\end{equation}
Using the fact that the following relations hold for $\tau$ and $\chi$ sufficiently small
\begin{equation}
    t_k\approx t, \quad t_{k+1}\approx t+\tau, \quad x_i\approx x, \quad x_{i\pm 1}\approx x\pm \chi,
    \label{eq:approx1}
\end{equation}
\begin{equation}
    n_{i}^{k}\approx n(x,t), \quad n_{i}^{k+1}\approx n(x,t+\tau), \quad c_{i}^{k}\approx c(x,t),
\end{equation}
\begin{equation}
   \quad \rho_n^k\approx \rho_n(t):=\int_\Omega n(x,t)\,\mathrm{d}x,  \quad K_{i}^{k}\approx K(x,t):=\int_\Omega \eta(x,x^\prime;\theta)c(x^\prime,t)\,\mathrm{d}x^\prime,
   \label{eq:approx2}
\end{equation}
where the function $\eta$ is defined via \eqref{eq:function_eta}, equation~\eqref{princ_mass} can be formally rewritten in the approximate form 
\begin{equation}
n(x,t+\tau)-n(x,t)=\tau \left(\alpha_n-\zeta_n K(x,t)-\mu_n\rho_n(t) \right)n(x,t).
\label{approx_form}
\end{equation}
If, in addiction, the function $n(x,t)$ is continuously differentiable with respect to the variable $t$, starting from equation~\eqref{approx_form}, and letting the time-step $\tau \to 0$, one formally obtains the following IDE for the tumour cell density $n(x,t)$:
$$ \partial_t n(x,t) =\alpha_n n(x,t)-\mu_n\rho_n(t)\, n(x,t)-\zeta_nK(x,t) n(x,t)\quad (x,t)\in  \Omega\times(0,t_f].$$ 
\subsection{Formal derivation of the PDE for the density of T cells $c(x,t)$}
When T  cell dynamics are governed by the rules described in Section~\ref{Dynamics of T  cells}, considering $i \in [1,\mathcal{N}-1]$, between time-steps $k$ and $k+1$ the principle of mass balance gives the following difference equation for the T cell density $c_{i}^{k}$:
\begin{equation}
\begin{aligned}
 c_i^{k+1}&=c_i^k(1-\tau\,\mu_c\rho_c^k)+\tau\,\alpha_cr^k_i \\
&+ \frac{\lambda}{2}\psi(w_i^k)\Big(c_{i+1}^k+c_{i-1}^k\Big) \\
&-\frac{\lambda}{2}\Big(\psi(w_{i-1}^k)+\psi(w_{i+1}^k)\Big)c_i^k\\
&+\frac{\nu}{2\phi_{\max}}\psi(w_i^k)\Big[\Big(\phi_i^k-\phi_{i-1}^k\Big)_{+}c_{i-1}^k\Big]\\
&+\frac{\nu}{2\phi_{\max}}\psi(w_i^k)\Big[\Big(\phi_i^k-\phi_{i+1}^k)\Big)_{+}c_{i+1}^k\Big] \\
&-\frac{\nu}{2\phi_{\max}}\psi(w_{i+1}^k)\Big[\Big(\phi_{i+1}^k-\phi_i^k\Big)_{+}c_i^k\Big]\\
&-\frac{\nu}{2\phi_{\max}}\psi(w_{i-1}^k)\Big[\Big(\phi_{i-1}^k-\phi_i^k\Big)_{+}c_i^k\Big].
\end{aligned}
\label{master2}
\end{equation}
Using the fact that relations \eqref{eq:approx1}-\eqref{eq:approx2} and the following relations 
\begin{equation*}
    c_i^k\approx c(x,t), \quad c_{i\pm 1}^k\approx c(x\pm\chi), \quad \rho_c^k\approx \rho_c(t):=\int_\Omega c(x,t)\,\mathrm{d}x,
    \label{approx_3}
\end{equation*}
\begin{equation*}
    \phi_i^k\approx \phi(x,t), \quad \phi_i^{k+1}\approx \phi(x,t+\tau), \quad \phi_{i\pm 1}^k\approx\phi(x\pm\chi),
    \label{approx_4}
\end{equation*}
\begin{equation*}
    w_i^k\approx w(x,t), \text{ with } w(x,t)=n(x,t)+c(x,t), \quad w_{i\pm 1}^k\approx w(x\pm\chi),
    \label{approx_4bis}
\end{equation*}
\begin{equation*}
    r_i^k\approx r(x,t):=\phi_{tot}(t)\mathbb{1}_\omega(x), \text{ with } \phi_{tot}(t):=\displaystyle\int_\Omega \phi(x,t)\,\mathrm{d}x
     \label{approx_5}
\end{equation*}
hold for $\tau$ and $\chi$ sufficiently small, equation~\eqref{master2} can be formally rewritten in the approximate form
\begin{equation*}
\begin{aligned}
c(x,t+\tau)&=c(x,t)(1-\tau\,\mu_c\rho_c(t))+\tau\,\alpha_cr(x,t) \\
&+ \frac{\lambda}{2}\psi(w(x,t))\Big(c(x+\chi,t)+c(x-\chi,t)\Big)\\
&-\frac{\lambda}{2}\Big(\psi(w(x-\chi,t))+\psi(w(x+\chi,t))\Big)c(x,t) \\
&+\frac{\nu}{2\phi_{\max}}\psi(w(x,t))\Big[\Big(\phi(x,t)-\phi(x-\chi,t)\Big)_{+}c(x-\chi,t)\Big]\\
&+\frac{\nu}{2\phi_{\max}}\psi(w(x,t))\Big[\Big(\phi(x,t)-\phi(x+\chi,t)\Big)_{+}c(x+\chi,t)\Big] \\
&-\frac{\nu}{2\phi_{\max}}\psi(w(x+\chi,t))\Big[\Big(\phi(x+\chi,t)-\phi(x,t)\Big)_{+}c(x,t)\Big]\\
&-\frac{\nu}{2\phi_{\max}}\psi(w(x-\chi,t))\Big[\Big(\phi(x-\chi,t)-\phi(x,t)\Big)_{+}c(x,t)\Big].
\end{aligned}
\label{master3}
\end{equation*}
Building on the methods employed in~\citep{bubba2020discrete}, letting $\tau \rightarrow 0$ and
$\chi \rightarrow 0$ in such a way that
\begin{equation*}
    \frac{\lambda}{2}\frac{\chi^2}{\tau}\rightarrow \beta_c \in \mathbb{R}_*^+ \quad \text{and} \quad  \frac{\nu}{2\phi_{\max}}\frac{\chi^2}{\tau}\rightarrow \gamma_c \in \mathbb{R}_*^+ \quad \text{as } \tau \rightarrow 0, \, \chi \rightarrow 0,
    \label{app:conditionsContMod}
\end{equation*}
after a little algebra, considering $(x,t)\in  \Omega\setminus\partial \Omega\times(0,t_f]$, we find
$$
   \partial_t c- \partial_x\Big[\beta_c \psi(w) \partial_x c-\gamma_c \psi(w)c\partial_x \phi -\beta_c c\psi^\prime(w)\partial_x w\Big]=  -\mu_c\rho_c(t)c+\alpha_cr  
  $$
where $\psi$ is given by \eqref{sqeezing_prob} and $w:=n+c$.
Moreover, zero-flux boundary conditions 
\begin{comment}
\begin{equation*}
\partial_x c(x,t)=0 \quad (x,t)\in \partial \Omega\times\mathbb{R}^*_+
\label{ap:boundaryCondition}
\end{equation*}
\end{comment}
easily follow from the fact that T-cell moves that require moving out of the spatial domain are not allowed. 
\subsection{Formal derivation of the balance equation for the chemoattractant concentration $\phi(x,t)$}
The formal derivation of the balance equation for the chemoattractant concentration $\phi(x,t)$ is obtained using the methods employed in~\citep{bubba2020discrete}.
\begin{comment}
In the case where the chemoattractant dynamics are governed by the rules described in Section~\ref{Dynamic of for the chemoattractant}, considering $(i,k) \in [1,I-1]\times\mathbb{N_0}$ and assuming that relations \eqref{eq:approx1}-\eqref{eq:approx2} and relations \eqref{approx_3}-\eqref{approx_4} hold for $\tau, \chi$ sufficiently small, the difference Equation~\eqref{chemoattr} can be formally rewritten in the approximate form
\begin{equation}
   \frac{\phi(x, t+\tau)-\phi(x,t)}{\tau}=\beta_{\phi}\frac{\phi(x-\chi,t)+\phi(x+\chi,t)-2\phi(x,t)}{\chi^2}+\alpha_{\phi}n(x,t)-\kappa_{\phi}\phi(x,t)
\end{equation}
If the function $\phi$ is continuously differentiable with respect to the variable $t$ and twice differentiable with respect to the variable $x$, letting $\tau, \chi \rightarrow 0$ in the above equation gives
\begin{equation}
    \partial_t \phi(x,t)-\beta_{\phi}\partial_{xx} \phi(x,t)=\alpha_{\phi}n(x,t)-\kappa_{\phi}\phi(x,t)
    \label{cont_chemo}
\end{equation}
which is the balance equation for the chemoattractant concentration $\phi$ posed on $\Omega\times \mathbb{R}^*_+$.
We complement Equation~\eqref{cont_chemo} with zero flux boundary conditions.
\end{comment}
\section{Details of numerical simulations}
The numerical simulations of our hybrid and continuum models are carried out on a two-dimensional domain and are performed in \textsc{Matlab}. 
\subsection{Details of numerical simulations of the hybrid model}
\label{app:sim1Ddiscrete}
\begin{figure}
\centering
\includegraphics[scale=0.7]{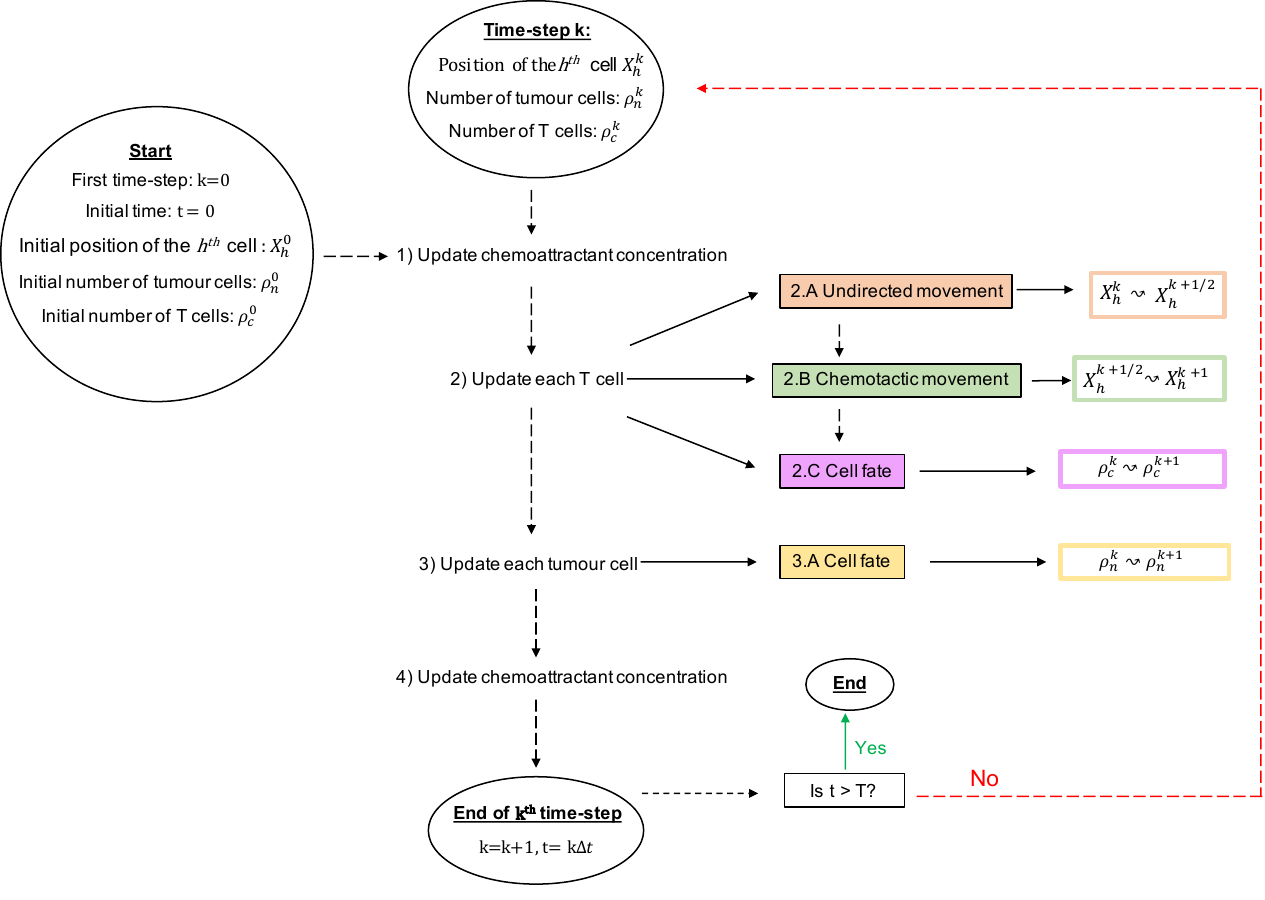}
\caption{Flowchart illustrating the computational procedure to simulate the hybrid model
in one-dimensional settings. A detailed summary of steps 2) and 3) is provided by the flowchart in Fig.~\ref{fig:schemaIBM2}. A similar procedure is used in two-dimensional settings.}
\label{fig:schemaIBM1}
\end{figure}

\begin{figure}
\centering
\includegraphics[scale=0.55,angle=-90]{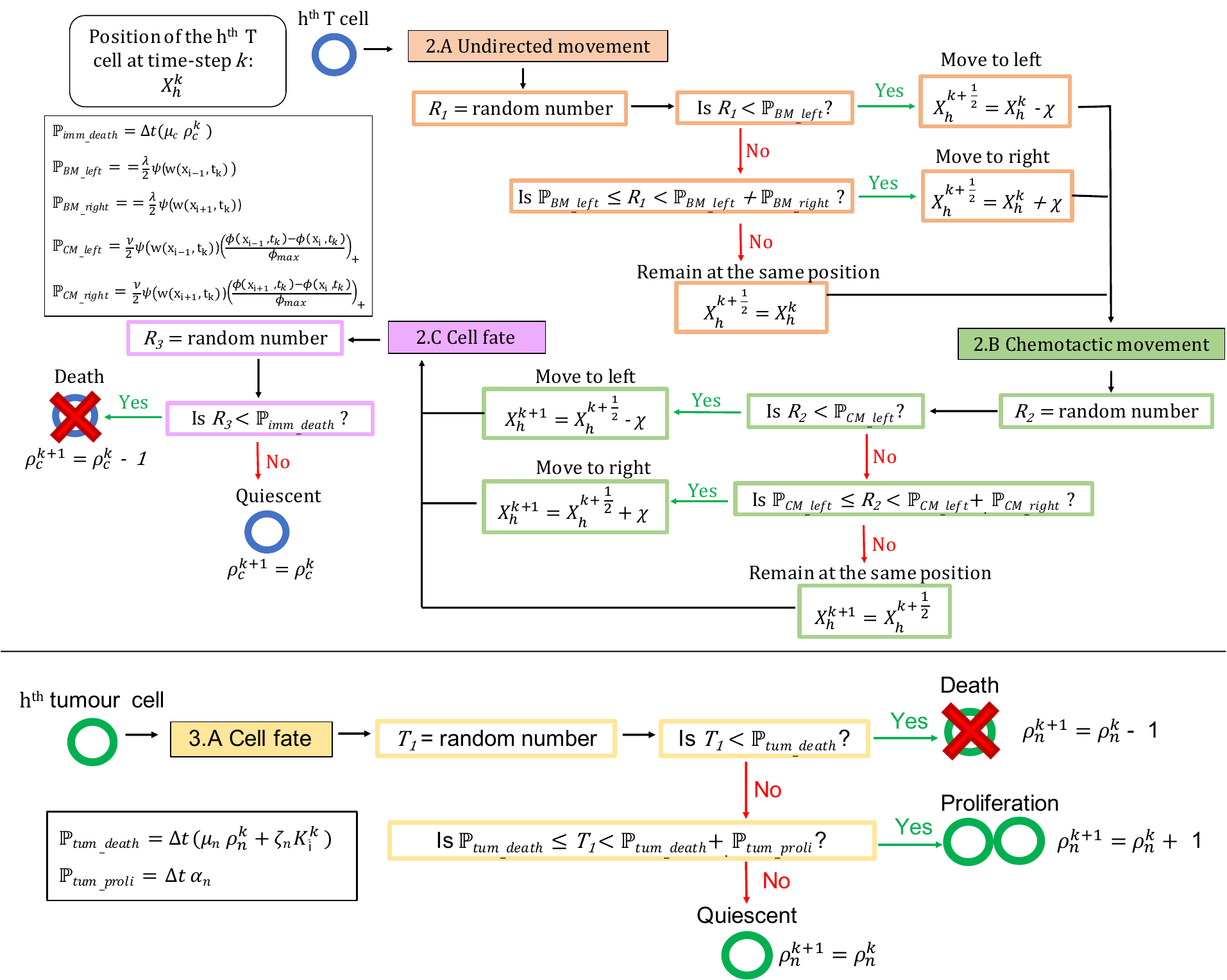}
\caption{Flowchart illustrating the detailed computational procedure followed to update the positions of every T cell, as well as the fate of each tumour cell and T cell
in one-dimensional settings. Analogous strategies are used in two-dimensional domains.}
\label{fig:schemaIBM2}
\end{figure}

The flowchart in Fig.~\ref{fig:schemaIBM1} illustrates the general computational procedure to carry out simulations of the hybrid model in one-dimensional settings, while the flowchart in Fig.~\ref{fig:schemaIBM2} provides further details of the computational procedure to simulate cell dynamics in one-dimensional settings. Analogous strategies are used in two-dimensional settings.
All random numbers mentioned in Fig. \ref{fig:schemaIBM2} are real numbers drawn from the standard uniform distribution on the interval $(0,1)$, which in our case are obtained using the  built-in \textsc{Matlab} function \textsc{rand}. %At each time-step $\tau$, the positions of the single T cells are updated following a procedure analogous to that explaned in Fig.~\ref{fig:schemaIBM1} and Fig.~\ref{fig:schemaIBM2}), with the only difference being that the T cells are allowed to move up and down as well.  

As summarised by Fig.~\ref{fig:schemaIBM2}, at any time-step, each T cell undergoes a three-phase process: Phase A) undirected, random movement according to the probabilities defined via \eqref{diff_left} and \eqref{diff_right}; Phase B) chemotaxis according to the probabilities defined via \eqref{chem_left} and \eqref{chem_right}; Phase C) death according to the probabilities defined via~\eqref{eq:inflow} and \eqref{eq:deathTcells}. We let then each tumour cell proliferate with the probability defined via \eqref{eq:tumProl}, die due to intra-tumour competition with the probability defined via \eqref{eq:ProbmortTum}, or die due to immune action with the probability defined via \eqref{eq:mortCompetition}. Finally, the tumour cell density
at every lattice site is computed via \eqref{eq:tumourDensity} and inserted into \eqref{chemoattr} in order to update the concentration of
the chemoattractant. 

In a two-dimensional setting, the positions of the single T cells are updated following a procedure analogous to that illustrated in Figs.~\ref{fig:schemaIBM1} and~\ref{fig:schemaIBM2}, with the only differences being that: T cells are allowed to move up and down as well; the concentration of the chemoattractant is updated through the two-dimensional analogue of \eqref{chemoattr}, where the operator $\mathcal{L}$ is defined as the finite-difference Laplacian on a two-dimensional regular lattice of step $\chi$; the tumour and T  cell densities are respectively computed via \eqref{eq:tumourDensity} and \eqref{eq:tcellDensity}.

\subsection{Details of numerical simulations of the continuum model}
\label{app:sim1Dcontinuum}
To construct numerical solutions of the IDE-PDE-PDE system \eqref{continuum_model}, we use a uniform discretisation consisting of $N^2 = 3721$ points of the square 
 $\Omega:= [0, 1]^2$ as
the computational domain of the independent variable $\textbf{x}\equiv  (x, y)$ (\textit{i.e.} $(x_{i}, y_{j}) = (i\Delta x, j\Delta x)$ with
$\Delta x=0.016$ and $i, j = 0,\dots,N)$ Moreover, we choose the time step $\Delta t= 10^{-4}$ and, unless stated otherwise, we perform numerical simulations for $15\times 10^4$ time-steps (\textit{i.e.} the final time of simulations is $t_f=15$).\\
The method for constructing numerical solutions of the IDE-PDE-PDE system \eqref{continuum_model} is based on a finite difference
scheme whereby the discretised dependent variables are
$$n_{i,j}^{k}:=n(x_{i},y_{j},t_k), \quad c_{i,j}^{k}:=c(x_{i},y_{j},t_k) \quad \text{and} \quad \phi_{i,j}^{k}:=\phi(x_{i},y_{j},t_k) . $$
We solve numerically the IDE~\eqref{continuum_model}$_1$ for $n$ and the PDE~\eqref{continuum_model}$_3$ for $\phi$ using the following schemes
$$ 
   \frac{n_{i,j}^{k+1}-n_{i,j}^{k}}{\Delta t}=\left(\alpha_n-\mu_n\rho_n^k-\zeta_nK_{i,j}^{k}\right)n_{i,j}^{k} \quad i,j=0,\dots,N,
$$
and
\begin{equation*}
  \begin{aligned}
\frac{\phi_{i,j}^{k+1}-\phi_{i,j}^{k}}{\Delta t}=&\beta_{\phi}\frac{\phi_{i-1,j}^k+\phi_{i+1,j}^k-2\phi_{i,j}^{k}}{(\Delta x)^2}+\beta_{\phi}\frac{\phi_{i,j-1}^k+\phi_{i,j+1}^k-2\phi_{i,j}^{k}}{(\Delta x)^2}\\
   &+\alpha_{\phi}n^k_{i,j}-\kappa_{\phi}\phi_{i,j}^{k}, \quad i,j=1,\dots,N-1,
   \end{aligned}
\end{equation*}
and impose zero-flux boundary conditions for $\phi$ by letting
$$\phi_{0,j}^{k+1}=\phi_{1,j}^{k+1} \quad \text{and} \quad \phi_{N,j}^{k+1}=\phi_{N-1,j}^{k+1}, \quad j=0,\dots, N $$
$$\phi_{i,0}^{k+1}=\phi_{i,1}^{k+1} \quad \text{and} \quad \phi_{i,N}^{k+1}=\phi_{i,N-1}^{k+1}, \quad i=0,\dots, N $$
Moreover, we solve numerically the PDE~\eqref{continuum_model}$_2$  for $c$ using the following explicit scheme, which is the same as the one used in~\citep{bubba2020discrete}, 
$$   
   \frac{c_{i,j}^{k+1}-c_{i,j}^{k}}{\Delta t}-\frac{F^{k}_{i+\frac{1}{2},j}+F^{k}_{i-\frac{1}{2},j}}{\Delta x}-\frac{F^{k}_{i,j+\frac{1}{2}}+F^{k}_{i,j-\frac{1}{2}}}{\Delta x}=r(\phi^k_{i,j})-\mu_c\rho_c^kc_{i,j}^{k}
$$
for  $i,j=0,\dots, N$, where
\begin{equation*}
\begin{aligned}
&F^{k}_{i+\frac{1}{2},j}:=\beta_c \psi(w^k_{i+\frac{1}{2},j}) \frac{c^k_{i+1,j}-c_{i,j}^{k}}{\Delta x}-\beta_c c^k_{i+\frac{1}{2},j}\psi^\prime(w^k_{i+\frac{1}{2},j})\frac{w^k_{i+1,j}-w_{i,j}^{k}}{\Delta x}\\
&-b^{k,+}_{i+\frac{1}{2},j}c_{i,j}^k\psi(w_{i+1,j}^k)+b^{k,-}_{i+\frac{1}{2},j}c_{i+1,j}^{k}\psi(w_{i,j}^k), \quad i=0, \dots N-1, \; j=0, \dots N,
\end{aligned}
\end{equation*}
\begin{equation*}
\begin{aligned}
&F^{k}_{i,j+\frac{1}{2}}:=\beta_c \psi(w^k_{i,j+\frac{1}{2}}) \frac{c^k_{i,j+1}-c_{i,j}^{k}}{\Delta x}-\beta_c c^k_{i,j+\frac{1}{2}}\psi^\prime(w^k_{i,j+\frac{1}{2}})\frac{w^k_{i,j+1}-w_{i,j}^{k}}{\Delta x}\\
&-b^{k,+}_{i,j+\frac{1}{2}}c_{i,j}^k\psi(w_{i,j+1}^k)+b^{k,-}_{i,j+\frac{1}{2}}c_{i,j+1}^{k+1}\psi(w_{i,j}^k), \qquad i=0, \dots N, \;  j=0, \dots N-1,
\end{aligned}
\end{equation*}
with 
$$w^k_{i+\frac{1}{2},j}:=\frac{w^k_{i+1,j}+w^k_{i,j}}{2},  \qquad w^k_{i,j+\frac{1}{2}}:=\frac{w^k_{i,j+1}+w^k_{i,j}}{2},$$
$$c^k_{i+\frac{1}{2},j}:=\frac{c^k_{i+1,j}+c^k_{i,j}}{2},  \qquad c^k_{i,j+\frac{1}{2}}:=\frac{c^k_{i,j+1}+c^k_{i,j}}{2},$$
$$b^{k}_{i+\frac{1}{2},j}:=\gamma_c\frac{\phi^k_{i+1,j}-\phi^k_{i,j}}{\Delta x}, \; b^{k,+}_{i+\frac{1}{2},j}=\max\left(0,b^{k}_{i+\frac{1}{2},j}\right), \; b^{k,-}_{i+\frac{1}{2},j}=\max\left(0,-b^{k}_{i+\frac{1}{2},j}\right)$$
and 
$$b^{k}_{i,j+\frac{1}{2}}:=\gamma_c\frac{\phi^k_{i,j+1}-\phi^k_{i,j}}{\Delta x}, \; b^{k,+}_{i,j+\frac{1}{2}}=\max\left(0,b^{k}_{i,j+\frac{1}{2}}\right), \; b^{k,-}_{i,j+\frac{1}{2}}=\max\left(0,-b^{k}_{i,j+\frac{1}{2}}\right)$$
The discrete fluxes $F^k_{i-\frac{1}{2},j}$ for $i=1,\dots, N,\, j=0,\dots, N$ and $F^k_{i,j-\frac{1}{2}}$ for $i=0,\dots, N, \,j=1,\dots, N$ are defined in an analogous way, and we impose zero-flux boundary
conditions by using the definitions
$$F^k_{0-\frac{1}{2},j}:=0 \quad \text{and} \quad  F^k_{N+\frac{1}{2},j}:=0, \quad \text{for} \quad j=0,\dots,N,$$
$$F^k_{i,0-\frac{1}{2}}:=0 \quad \text{and} \quad  F^k_{i,N+\frac{1}{2}}:=0, \quad \text{for} \quad i=0,\dots,N.$$
	 \end{appendices}

\bibliography{sn-bibliography}% common bib file

\begin{thebibliography}{74}
\providecommand{\natexlab}[1]{#1}
\providecommand{\url}[1]{\texttt{#1}}
\expandafter\ifx\csname urlstyle\endcsname\relax
  \providecommand{\doi}[1]{doi: #1}\else
  \providecommand{\doi}{doi: \begingroup \urlstyle{rm}\Url}\fi

\bibitem[Aguad{\'e}-Gorgori{\'o} and Sol{\'e}(2020)]{aguade2020tumour}
G.~Aguad{\'e}-Gorgori{\'o} and R.~Sol{\'e}.
\newblock Tumour neoantigen heterogeneity thresholds provide a time window for
  combination immunotherapy.
\newblock \emph{J. R. Soc. Interface}, 17\penalty0 (171):\penalty0 20200736,
  2020.

\bibitem[Al-Tameemi et~al.(2012)Al-Tameemi, Chaplain, and
  d'Onofrio]{al2012evasion}
M.~Al-Tameemi, M.~Chaplain, and A.~d'Onofrio.
\newblock Evasion of tumours from the control of the immune system:
  consequences of brief encounters.
\newblock \emph{Biol. Direct}, 7\penalty0 (1):\penalty0 31, 2012.

\bibitem[Almeida et~al.(2021)Almeida, Audebert, Leschiera, and
  Lorenzi]{almeida2021discrete}
L.~Almeida, C.~Audebert, E.~Leschiera, and T.~Lorenzi.
\newblock Discrete and continuum models for the coevolutionary dynamics between
  {CD8+ cytotoxic T} lymphocytes and tumour cells.
\newblock \emph{arXiv preprint arXiv:2109.09568}, 2021.

\bibitem[Almuallem et~al.(2021)Almuallem, Trucu, and
  Eftimie]{almuallem2021oncolytic}
N.~Almuallem, D.~Trucu, and R.~Eftimie.
\newblock Oncolytic viral therapies and the delicate balance between
  virus-macrophage-tumour interactions: A mathematical approach.
\newblock \emph{Math. Biosci. Eng.}, 18\penalty0 (1):\penalty0 764--799, 2021.

\bibitem[Angell and Galon(2013)]{angell2013immune}
H.~Angell and J.~Galon.
\newblock From the immune contexture to the immunoscore: the role of prognostic
  and predictive immune markers in cancer.
\newblock \emph{Curr. Opin. Immunol.}, 25\penalty0 (2):\penalty0 261--267,
  2013.

\bibitem[Atsou et~al.(2020)Atsou, Anju{\`e}re, Braud, and
  Goudon]{atsou2020size}
K.~Atsou, F.~Anju{\`e}re, V.~M. Braud, and T.~Goudon.
\newblock A size and space structured model describing interactions of tumor
  cells with immune cells reveals cancer persistent equilibrium states in
  tumorigenesis.
\newblock \emph{J Theor Biol}, 490:\penalty0 110163, 2020.

\bibitem[Basu et~al.(2016)Basu, Whitlock, Husson, Le~Floc'h, Jin, Oyler-Yaniv,
  Dotiwala, Giannone, Hivroz, Biais, et~al.]{basu2016cytotoxic}
R.~Basu, B.~M. Whitlock, J.~Husson, A.~Le~Floc'h, W.~Jin, A.~Oyler-Yaniv,
  F.~Dotiwala, G.~Giannone, C.~Hivroz, N.~Biais, et~al.
\newblock Cytotoxic {T cells} use mechanical force to potentiate target cell
  killing.
\newblock \emph{Cell}, 165\penalty0 (1):\penalty0 100--110, 2016.

\bibitem[Boissonnas et~al.(2007)Boissonnas, Fetler, Zeelenberg, Hugues, and
  Amigorena]{boissonnas2007vivo}
A.~Boissonnas, L.~Fetler, I.~S. Zeelenberg, S.~Hugues, and S.~Amigorena.
\newblock In vivo imaging of cytotoxic {T} cell infiltration and elimination of
  a solid tumor.
\newblock \emph{J. Exp. Med.}, 204\penalty0 (2):\penalty0 345--356, 2007.

\bibitem[Bubba et~al.(2020)Bubba, Lorenzi, and Macfarlane]{bubba2020discrete}
F.~Bubba, T.~Lorenzi, and F.~R. Macfarlane.
\newblock From a discrete model of chemotaxis with volume-filling to a
  generalized {Patlak--Keller--Segel} model.
\newblock \emph{Proc R Soc Lond A}, 476\penalty0 (2237):\penalty0 20190871,
  2020.

\bibitem[Byrne and Drasdo(2009)]{byrne2009individual}
H.~Byrne and D.~Drasdo.
\newblock Individual-based and continuum models of growing cell populations: a
  comparison.
\newblock \emph{J. Math. Biol.}, 58\penalty0 (4):\penalty0 657--687, 2009.

\bibitem[Cattani et~al.(2010)Cattani, Ciancio, and d'Onofrio]{Cattani2010}
C.~Cattani, A.~Ciancio, and A.~d'Onofrio.
\newblock Metamodeling the learning--hiding competition between tumours and the
  immune system: a kinematic approach.
\newblock \emph{Math. Comput. Model. Dyn. Syst.}, 52\penalty0 (1):\penalty0
  62--69, 2010.

\bibitem[Champagnat et~al.(2008)Champagnat, Ferri{\`e}re, and
  M{\'e}l{\'e}ard]{champagnat2008individual}
N.~Champagnat, R.~Ferri{\`e}re, and S.~M{\'e}l{\'e}ard.
\newblock From individual stochastic processes to macroscopic models in
  adaptive evolution.
\newblock \emph{Stoch. Models}, 24\penalty0 (sup1):\penalty0 2--44, 2008.

\bibitem[Chisholm et~al.(2016)Chisholm, Lorenzi, Desvillettes, and
  Hughes]{chisholm2016evolutionary}
R.~H. Chisholm, T.~Lorenzi, L.~Desvillettes, and B.~D. Hughes.
\newblock Evolutionary dynamics of phenotype-structured populations: from
  individual-level mechanisms to population-level consequences.
\newblock \emph{Z. Angew. Math. Phys.}, 67\penalty0 (4):\penalty0 100, 2016.

\bibitem[Christophe et~al.(2015)Christophe, M{\"u}ller, Rodrigues, Petit,
  Cattiaux, Dupr{\'e}, Gadat, and Valitutti]{christophe2015biased}
C.~Christophe, S.~M{\"u}ller, M.~Rodrigues, A.-E. Petit, P.~Cattiaux,
  L.~Dupr{\'e}, S.~Gadat, and S.~Valitutti.
\newblock A biased competition theory of cytotoxic {T lymphocyte interaction
  with tumor nodules}.
\newblock \emph{PloS ONE}, 10\penalty0 (3), 2015.

\bibitem[Cooper and Kim(2014)]{cooper}
A.~K. Cooper and P.~S. Kim.
\newblock A cellular automata and a partial differential equation model of
  tumor--immune dynamics and chemotaxis.
\newblock In A.~Eladdadi, P.~Kim, and D.~Mallet, editors, \emph{Mathematical
  Models of Tumor-Immune System Dynamics}, pages 21--46, New York, NY, 2014.
  Springer New York.

\bibitem[Coulie et~al.(2014)Coulie, Van~den Eynde, Van Der~Bruggen, and
  Boon]{coulie2014tumour}
P.~G. Coulie, B.~J. Van~den Eynde, P.~Van Der~Bruggen, and T.~Boon.
\newblock Tumour antigens recognized by {T} lymphocytes: at the core of cancer
  immunotherapy.
\newblock \emph{Nat. Rev. Cancer}, 14\penalty0 (2):\penalty0 135--146, 2014.

\bibitem[Delitala and Lorenzi(2013)]{delitala2013recognition}
M.~Delitala and T.~Lorenzi.
\newblock Recognition and learning in a mathematical model for immune response
  against cancer.
\newblock \emph{Discrete Contin. Dyn. Syst. - B}, 18\penalty0 (4), 2013.

\bibitem[Eftimie et~al.(2011)Eftimie, Bramson, and
  Earn]{eftimie2011interactions}
R.~Eftimie, J.~L. Bramson, and D.~J. Earn.
\newblock Interactions between the immune system and cancer: a brief review of
  non-spatial mathematical models.
\newblock \emph{Bull. Math. Biol.}, 73\penalty0 (1):\penalty0 2--32, 2011.

\bibitem[Galon and Bruni(2019)]{galon2019approaches}
J.~Galon and D.~Bruni.
\newblock Approaches to treat immune hot, altered and cold tumours with
  combination immunotherapies.
\newblock \emph{Nat. Rev. Drug Discov.}, 18\penalty0 (3):\penalty0 197--218,
  2019.

\bibitem[Galon et~al.(2006)Galon, Costes, Sanchez-Cabo, Kirilovsky, Mlecnik,
  Lagorce-Pag{\`e}s, Tosolini, Camus, Berger, Wind, et~al.]{galon2006type}
J.~Galon, A.~Costes, F.~Sanchez-Cabo, A.~Kirilovsky, B.~Mlecnik,
  C.~Lagorce-Pag{\`e}s, M.~Tosolini, M.~Camus, A.~Berger, P.~Wind, et~al.
\newblock Type, density, and location of immune cells within human colorectal
  tumors predict clinical outcome.
\newblock \emph{Science}, 313\penalty0 (5795):\penalty0 1960--1964, 2006.

\bibitem[Galon et~al.(2016)Galon, Fox, Bifulco, Masucci, Rau, Botti, Marincola,
  Ciliberto, Pages, Ascierto, et~al.]{galon2016immunoscore}
J.~Galon, B.~Fox, C.~Bifulco, G.~Masucci, T.~Rau, G.~Botti, F.~Marincola,
  G.~Ciliberto, F.~Pages, P.~Ascierto, et~al.
\newblock Immunoscore and {Immunoprofiling in cancer: an update from the
  melanoma and immunotherapy bridge 2015}, 2016.

\bibitem[Gandhi et~al.(2018)Gandhi, Rodr{\'\i}guez-Abreu, Gadgeel, Esteban,
  Felip, De~Angelis, Domine, Clingan, Hochmair, Powell,
  et~al.]{gandhi2018pembrolizumab}
L.~Gandhi, D.~Rodr{\'\i}guez-Abreu, S.~Gadgeel, E.~Esteban, E.~Felip,
  F.~De~Angelis, M.~Domine, P.~Clingan, M.~J. Hochmair, S.~F. Powell, et~al.
\newblock Pembrolizumab plus chemotherapy in metastatic non--small-cell lung
  cancer.
\newblock \emph{N. Engl. J. Med.}, 378\penalty0 (22):\penalty0 2078--2092,
  2018.

\bibitem[Gong et~al.(2017)Gong, Milberg, Wang, Vicini, Narwal, Roskos, and
  Popel]{gong2017computational}
C.~Gong, O.~Milberg, B.~Wang, P.~Vicini, R.~Narwal, L.~Roskos, and A.~S. Popel.
\newblock A computational multiscale agent-based model for simulating
  spatio-temporal tumour immune response to {PD1 and PDL1} inhibition.
\newblock \emph{J. R. Soc. Interface}, 14\penalty0 (134):\penalty0 20170320,
  2017.

\bibitem[Gorbachev et~al.(2007)Gorbachev, Kobayashi, Kudo, Tannenbaum, Finke,
  Shu, Farber, and Fairchild]{gorbachev2007cxc}
A.~V. Gorbachev, H.~Kobayashi, D.~Kudo, C.~S. Tannenbaum, J.~H. Finke, S.~Shu,
  J.~M. Farber, and R.~L. Fairchild.
\newblock Cxc chemokine ligand 9/monokine induced by ifn-$\gamma$ production by
  tumor cells is critical for t cell-mediated suppression of cutaneous tumors.
\newblock \emph{J Immunol}, 178\penalty0 (4):\penalty0 2278--2286, 2007.

\bibitem[Griffiths et~al.(2020)Griffiths, Wallet, Pflieger, Stenehjem, Liu,
  Cosgrove, Leggett, McQuerry, Shrestha, Rossetti,
  et~al.]{griffiths2020circulating}
J.~I. Griffiths, P.~Wallet, L.~T. Pflieger, D.~Stenehjem, X.~Liu, P.~A.
  Cosgrove, N.~A. Leggett, J.~A. McQuerry, G.~Shrestha, M.~Rossetti, et~al.
\newblock Circulating immune cell phenotype dynamics reflect the strength of
  tumor--immune cell interactions in patients during immunotherapy.
\newblock \emph{Proc. Natl. Acad. Sci. U.S.A.}, 117\penalty0 (27):\penalty0
  16072--16082, 2020.

\bibitem[Halle et~al.(2016)Halle, Keyser, Stahl, Busche, Marquardt, Zheng,
  Galla, Heissmeyer, Heller, Boelter, et~al.]{halle2016vivo}
S.~Halle, K.~A. Keyser, F.~R. Stahl, A.~Busche, A.~Marquardt, X.~Zheng,
  M.~Galla, V.~Heissmeyer, K.~Heller, J.~Boelter, et~al.
\newblock In vivo killing capacity of cytotoxic {T cells is limited and
  involves dynamic interactions and } cooperativity.
\newblock \emph{Immunity}, 44\penalty0 (2):\penalty0 233--245, 2016.

\bibitem[Handel et~al.(2020)Handel, La~Gruta, and Thomas]{handel2020simulation}
A.~Handel, N.~L. La~Gruta, and P.~G. Thomas.
\newblock Simulation modelling for immunologists.
\newblock \emph{Nat. Rev. Immunol.}, 20\penalty0 (3):\penalty0 186--195, 2020.

\bibitem[Hegde et~al.(2016)Hegde, Karanikas, and Evers]{hegde2016and}
P.~S. Hegde, V.~Karanikas, and S.~Evers.
\newblock The where, the when, and the how of immune monitoring for cancer
  immunotherapies in the era of checkpoint inhibition.
\newblock \emph{Clinical Cancer Res.}, 22\penalty0 (8):\penalty0 1865--1874,
  2016.

\bibitem[Hellmann et~al.(2018)Hellmann, Ciuleanu, Pluzanski, Lee, Otterson,
  Audigier-Valette, Minenza, Linardou, Burgers, Salman,
  et~al.]{hellmann2018nivolumab}
M.~D. Hellmann, T.-E. Ciuleanu, A.~Pluzanski, J.~S. Lee, G.~A. Otterson,
  C.~Audigier-Valette, E.~Minenza, H.~Linardou, S.~Burgers, P.~Salman, et~al.
\newblock Nivolumab plus ipilimumab in lung cancer with a high tumor mutational
  burden.
\newblock \emph{N. Engl. J. Med.}, 378\penalty0 (22):\penalty0 2093--2104,
  2018.

\bibitem[Huang et~al.(2017)Huang, Postow, Orlowski, Mick, Bengsch, Manne, Xu,
  Harmon, Giles, Wenz, et~al.]{huang2017t}
A.~C. Huang, M.~A. Postow, R.~J. Orlowski, R.~Mick, B.~Bengsch, S.~Manne,
  W.~Xu, S.~Harmon, J.~R. Giles, B.~Wenz, et~al.
\newblock T-cell invigoration to tumour burden ratio associated with
  anti-{PD-1} response.
\newblock \emph{Nature}, 545\penalty0 (7652):\penalty0 60--65, 2017.

\bibitem[Hughes et~al.(1995)]{hughes1995random}
B.~D. Hughes et~al.
\newblock \emph{Random walks and random environments: random walks}, volume~1.
\newblock Oxford University Press, 1995.

\bibitem[Iwai et~al.(2002)Iwai, Ishida, Tanaka, Okazaki, Honjo, and
  Minato]{iwai2002involvement}
Y.~Iwai, M.~Ishida, Y.~Tanaka, T.~Okazaki, T.~Honjo, and N.~Minato.
\newblock Involvement of {PD-L1} on tumor cells in the escape from host immune
  system and tumor immunotherapy by {PD-L1} blockade.
\newblock \emph{Proc. Natl. Acad. Sci. U.S.A.}, 99\penalty0 (19):\penalty0
  12293--12297, 2002.

\bibitem[Jarrett et~al.(2020)Jarrett, Faghihi, Hormuth, Lima, Virostko, Biros,
  Patt, and Yankeelov]{jarrett2020optimal}
A.~M. Jarrett, D.~Faghihi, D.~A. Hormuth, E.~A. Lima, J.~Virostko, G.~Biros,
  D.~Patt, and T.~E. Yankeelov.
\newblock Optimal control theory for personalized therapeutic regimens in
  oncology: Background, history, challenges, and opportunities.
\newblock \emph{J. Clin. Med.}, 9\penalty0 (5):\penalty0 1314, 2020.

\bibitem[Johnston et~al.(2015)Johnston, Simpson, and
  Baker]{johnston2015modelling}
S.~T. Johnston, M.~J. Simpson, and R.~E. Baker.
\newblock Modelling the movement of interacting cell populations: a moment
  dynamics approach.
\newblock \emph{J Theor Biol}, 370:\penalty0 81--92, 2015.

\bibitem[Kather et~al.(2017)Kather, Poleszczuk, Suarez-Carmona, Krisam,
  Charoentong, Valous, Weis, Tavernar, Leiss, Herpel, et~al.]{kather2017silico}
J.~N. Kather, J.~Poleszczuk, M.~Suarez-Carmona, J.~Krisam, P.~Charoentong,
  N.~A. Valous, C.-A. Weis, L.~Tavernar, F.~Leiss, E.~Herpel, et~al.
\newblock {\textit{In silico}} modeling of immunotherapy and stroma-targeting
  therapies in human colorectal cancer.
\newblock \emph{Cancer Res.}, 77\penalty0 (22):\penalty0 6442--6452, 2017.

\bibitem[Kato et~al.(2017)Kato, Yaguchi, Iwata, Morii, Nakagawa, Nishimura, and
  Kawakami]{kato2017prospects}
D.~Kato, T.~Yaguchi, T.~Iwata, K.~Morii, T.~Nakagawa, R.~Nishimura, and
  Y.~Kawakami.
\newblock Prospects for personalized combination immunotherapy for solid tumors
  based on adoptive cell therapies and immune checkpoint blockade therapies.
\newblock \emph{Nihon Rinsho Meneki Gakkai Kaishi}, 40\penalty0 (1):\penalty0
  68--77, 2017.

\bibitem[Kim and Lee(2012)]{kim2012modeling}
P.~S. Kim and P.~P. Lee.
\newblock Modeling protective anti-tumor immunity via preventative cancer
  vaccines using a hybrid agent-based and delay differential equation approach.
\newblock \emph{PLoS Comput. Biol.}, 8\penalty0 (10):\penalty0 e1002742, 2012.

\bibitem[Kolev(2003)]{Kolev2003}
M.~Kolev.
\newblock Mathematical modeling of the competition between acquired immunity
  and cancer.
\newblock \emph{Int. J. Appl. Math. Comput. Sci.}, 13:\penalty0 289--296, 2003.

\bibitem[Konstorum et~al.(2017)Konstorum, Vella, Adler, and
  Laubenbacher]{konstorum2017addressing}
A.~Konstorum, A.~T. Vella, A.~J. Adler, and R.~C. Laubenbacher.
\newblock Addressing current challenges in cancer immunotherapy with
  mathematical and computational modelling.
\newblock \emph{J. R. Soc. Interface}, 14\penalty0 (131):\penalty0 20170150,
  2017.

\bibitem[Kuznetsov and Knott(2001)]{Kuznetsov2001}
V.~A. Kuznetsov and G.~D. Knott.
\newblock {Modeling Tumor Regrowth and Immunotherapy}.
\newblock \emph{Math. Comput. Model.}, 33\penalty0 (12):\penalty0 1275--1287,
  2001.

\bibitem[Kuznetsov et~al.(1994)Kuznetsov, Makalkin, Taylor, and
  Perelson]{kuznetsov1994nonlinear}
V.~A. Kuznetsov, I.~A. Makalkin, M.~A. Taylor, and A.~S. Perelson.
\newblock Nonlinear dynamics of immunogenic tumors: parameter estimation and
  global bifurcation analysis.
\newblock \emph{Bull. Math. Biol.}, 56\penalty0 (2):\penalty0 295--321, 1994.

\bibitem[Leschiera et~al.(2022)Leschiera, Lorenzi, Shen, Almeida, and
  Audebert]{leschiera2022mathematical}
E.~Leschiera, T.~Lorenzi, S.~Shen, L.~Almeida, and C.~Audebert.
\newblock A mathematical model to study the impact of intra-tumour
  heterogeneity on anti-tumour {CD8+ T} cell immune response.
\newblock \emph{J Theor Biol}, page 111028, 2022.

\bibitem[Lin~Erickson et~al.(2009)Lin~Erickson, Wise, Fleming, Baird, Lateef,
  Molinaro, Teboh-Ewungkem, and de~Pillis]{LinErikson2009}
A.~H. Lin~Erickson, A.~Wise, S.~Fleming, M.~Baird, Z.~Lateef, A.~Molinaro,
  M.~Teboh-Ewungkem, and L.~G. de~Pillis.
\newblock A preliminary mathematical model of skin dendritic cell trafficking
  and induction of {T} cell immunity.
\newblock \emph{Discrete Contin. Dyn. Syst. - B}, 12:\penalty0 323--336, 2009.

\bibitem[Lorenzi(2022)]{lorenzi2022cancer}
T.~Lorenzi.
\newblock Cancer modelling as fertile ground for new mathematical challenges.
  comment on" improving cancer treatments via dynamical biophysical models" by
  m. kuznetsov, j. clairambault \& v. volpert.
\newblock \emph{Phys Life Rev}, 40:\penalty0 3--5, 2022.

\bibitem[Lorenzi et~al.(2015)Lorenzi, Chisholm, Melensi, Lorz, and
  Delitala]{lorenzi2015mathematical}
T.~Lorenzi, R.~H. Chisholm, M.~Melensi, A.~Lorz, and M.~Delitala.
\newblock Mathematical model reveals how regulating the three phases of
  {T-cell} response could counteract immune evasion.
\newblock \emph{Immunology}, 146\penalty0 (2):\penalty0 271--280, 2015.

\bibitem[{\L}uksza et~al.(2017){\L}uksza, Riaz, Makarov, Balachandran,
  Hellmann, Solovyov, Rizvi, Merghoub, Levine, Chan,
  et~al.]{luksza2017neoantigen}
M.~{\L}uksza, N.~Riaz, V.~Makarov, V.~P. Balachandran, M.~D. Hellmann,
  A.~Solovyov, N.~A. Rizvi, T.~Merghoub, A.~J. Levine, T.~A. Chan, et~al.
\newblock A neoantigen fitness model predicts tumour response to checkpoint
  blockade immunotherapy.
\newblock \emph{Nature}, 551\penalty0 (7681):\penalty0 517--520, 2017.

\bibitem[Macfarlane et~al.(2018)Macfarlane, Lorenzi, and
  Chaplain]{macfarlane2018modelling}
F.~R. Macfarlane, T.~Lorenzi, and M.~A. Chaplain.
\newblock Modelling the immune response to cancer: an individual-based approach
  accounting for the difference in movement between inactive and activated {T}
  cells.
\newblock \emph{Bull. Math. Biol.}, 80\penalty0 (6):\penalty0 1539--1562, 2018.

\bibitem[Macfarlane et~al.(2019)Macfarlane, Chaplain, and
  Lorenzi]{macfarlane2019stochastic}
F.~R. Macfarlane, M.~A. Chaplain, and T.~Lorenzi.
\newblock A stochastic individual-based model to explore the role of spatial
  interactions and antigen recognition in the immune response against solid
  tumours.
\newblock \emph{J Theor Biol}, 480:\penalty0 43--55, 2019.

\bibitem[Macfarlane et~al.(2020)Macfarlane, Chaplain, and
  Lorenzi]{macfarlane2020hybrid}
F.~R. Macfarlane, M.~A. Chaplain, and T.~Lorenzi.
\newblock A hybrid discrete-continuum approach to model turing pattern
  formation.
\newblock \emph{Math. Biosci. Eng.}, 17\penalty0 (6):\penalty0 7442--7479,
  2020.

\bibitem[Maini et~al.(1997)Maini, Painter, and Chau]{maini1997spatial}
P.~Maini, K.~Painter, and H.~P. Chau.
\newblock Spatial pattern formation in chemical and biological systems.
\newblock \emph{J. Chem. Soc. Faraday Trans.}, 93\penalty0 (20):\penalty0
  3601--3610, 1997.

\bibitem[Makaryan et~al.(2020)Makaryan, Cess, and Finley]{makaryan2020modeling}
S.~Z. Makaryan, C.~G. Cess, and S.~D. Finley.
\newblock Modeling immune cell behavior across scales in cancer.
\newblock \emph{Wiley Interdiscip. Rev. Syst. Biol. Med.}, 12\penalty0
  (4):\penalty0 e1484, 2020.

\bibitem[Mallet and De~Pillis(2006)]{mallet2006cellular}
D.~G. Mallet and L.~G. De~Pillis.
\newblock A cellular automata model of tumor--immune system interactions.
\newblock \emph{J Theor Biol}, 239\penalty0 (3):\penalty0 334--350, 2006.

\bibitem[MATLAB(2020)]{MATLAB:2021}
MATLAB.
\newblock \emph{9.9.0.1570001 (R2020b)}.
\newblock The MathWorks Inc., Natick, Massachusetts, 2020.

\bibitem[Matzavinos et~al.(2004)Matzavinos, Chaplain, and
  Kuznetsov]{matzavinos2004mathematical}
A.~Matzavinos, M.~A. Chaplain, and V.~A. Kuznetsov.
\newblock Mathematical modelling of the spatio-temporal response of cytotoxic
  {{T}-lymphocytes to a solid tumour}.
\newblock \emph{Math. Med. Biol.}, 21\penalty0 (1):\penalty0 1--34, 2004.

\bibitem[McGranahan et~al.(2016)McGranahan, Furness, Rosenthal, Ramskov,
  Lyngaa, Saini, Jamal-Hanjani, Wilson, Birkbak, Hiley,
  et~al.]{mcgranahan2016clonal}
N.~McGranahan, A.~J. Furness, R.~Rosenthal, S.~Ramskov, R.~Lyngaa, S.~K. Saini,
  M.~Jamal-Hanjani, G.~A. Wilson, N.~J. Birkbak, C.~T. Hiley, et~al.
\newblock Clonal neoantigens elicit {T} cell immunoreactivity and sensitivity
  to immune checkpoint blockade.
\newblock \emph{Science}, 351\penalty0 (6280):\penalty0 1463--1469, 2016.

\bibitem[Miller et~al.(2003)Miller, Wei, Cahalan, and
  Parker]{miller2003autonomous}
M.~J. Miller, S.~H. Wei, M.~D. Cahalan, and I.~Parker.
\newblock Autonomous {T cell} trafficking examined in vivo with intravital
  two-photon microscopy.
\newblock \emph{Proc. Natl. Acad. Sci. U.S.A.}, 100\penalty0 (5):\penalty0
  2604--2609, 2003.

\bibitem[Motzer et~al.(2018)Motzer, Tannir, McDermott, Frontera, Melichar,
  Choueiri, Plimack, Barth{\'e}l{\'e}my, Porta, George,
  et~al.]{motzer2018nivolumab}
R.~J. Motzer, N.~M. Tannir, D.~F. McDermott, O.~A. Frontera, B.~Melichar, T.~K.
  Choueiri, E.~R. Plimack, P.~Barth{\'e}l{\'e}my, C.~Porta, S.~George, et~al.
\newblock Nivolumab plus ipilimumab versus sunitinib in advanced renal-cell
  carcinoma.
\newblock \emph{N. Engl. J. Med.}, 2018.

\bibitem[Painter(2019)]{painter2019mathematical}
K.~J. Painter.
\newblock Mathematical models for chemotaxis and their applications in
  self-organisation phenomena.
\newblock \emph{Journal of theoretical biology}, 481:\penalty0 162--182, 2019.

\bibitem[Painter and Hillen(2002)]{painter2002volume}
K.~J. Painter and T.~Hillen.
\newblock Volume-filling and quorum-sensing in models for chemosensitive
  movement.
\newblock \emph{Can. Appl. Math. Quart}, 10\penalty0 (4):\penalty0 501--543,
  2002.

\bibitem[Pitt et~al.(2016)Pitt, Marabelle, Eggermont, Soria, Kroemer, and
  Zitvogel]{pitt2016targeting}
J.~Pitt, A.~Marabelle, A.~Eggermont, J.-C. Soria, G.~Kroemer, and L.~Zitvogel.
\newblock Targeting the tumor microenvironment: removing obstruction to
  anticancer immune responses and immunotherapy.
\newblock \emph{Ann. Oncol.}, 27\penalty0 (8):\penalty0 1482--1492, 2016.

\bibitem[Rabinovich et~al.(2007)Rabinovich, Gabrilovich, and
  Sotomayor]{rabinovich2007immunosuppressive}
G.~A. Rabinovich, D.~Gabrilovich, and E.~M. Sotomayor.
\newblock Immunosuppressive strategies that are mediated by tumor cells.
\newblock \emph{Annu. Rev. Immunol.}, 25:\penalty0 267--296, 2007.

\bibitem[Ribas and Wolchok(2018)]{ribas2018cancer}
A.~Ribas and J.~D. Wolchok.
\newblock Cancer immunotherapy using checkpoint blockade.
\newblock \emph{Science}, 359\penalty0 (6382):\penalty0 1350--1355, 2018.

\bibitem[Slaney et~al.(2014)Slaney, Kershaw, and Darcy]{slaney2014trafficking}
C.~Y. Slaney, M.~H. Kershaw, and P.~K. Darcy.
\newblock Trafficking of t cells into tumors.
\newblock \emph{Cancer Res.}, 74\penalty0 (24):\penalty0 7168--7174, 2014.

\bibitem[Spranger et~al.(2015)Spranger, Bao, and
  Gajewski]{spranger2015melanoma}
S.~Spranger, R.~Bao, and T.~F. Gajewski.
\newblock Melanoma-intrinsic $\beta$-catenin signalling prevents anti-tumour
  immunity.
\newblock \emph{Nature}, 523\penalty0 (7559):\penalty0 231--235, 2015.

\bibitem[Takayanagi and Ohuchi(2001)]{Takayanagi2001}
T.~Takayanagi and A.~Ohuchi.
\newblock A mathematical analysis of the interactions between immunogenic tumor
  cells and cytotoxic \uppercase{T} lymphocytes.
\newblock \emph{Microbiol. Immunol.}, 45\penalty0 (10):\penalty0 709--715,
  2001.

\bibitem[Tian et~al.(2017)Tian, Goldstein, Wang, Ching~Lo, Sun~Kim, Welte,
  Sheng, Dobrolecki, Zhang, Putluri, et~al.]{tian2017mutual}
L.~Tian, A.~Goldstein, H.~Wang, H.~Ching~Lo, I.~Sun~Kim, T.~Welte, K.~Sheng,
  L.~E. Dobrolecki, X.~Zhang, N.~Putluri, et~al.
\newblock Mutual regulation of tumour vessel normalization and
  immunostimulatory reprogramming.
\newblock \emph{Nature}, 544\penalty0 (7649):\penalty0 250--254, 2017.

\bibitem[Topalian et~al.(2012)Topalian, Hodi, Brahmer, Gettinger, Smith,
  McDermott, Powderly, Carvajal, Sosman, Atkins, et~al.]{topalian2012safety}
S.~L. Topalian, F.~S. Hodi, J.~R. Brahmer, S.~N. Gettinger, D.~C. Smith, D.~F.
  McDermott, J.~D. Powderly, R.~D. Carvajal, J.~A. Sosman, M.~B. Atkins, et~al.
\newblock Safety, activity, and immune correlates of {anti--PD-1} antibody in
  cancer.
\newblock \emph{N. Engl. J. Med.}, 366\penalty0 (26):\penalty0 2443--2454,
  2012.

\bibitem[Tumeh et~al.(2014)Tumeh, Harview, Yearley, Shintaku, Taylor, Robert,
  Chmielowski, Spasic, Henry, Ciobanu, et~al.]{tumeh2014pd}
P.~C. Tumeh, C.~L. Harview, J.~H. Yearley, I.~P. Shintaku, E.~J. Taylor,
  L.~Robert, B.~Chmielowski, M.~Spasic, G.~Henry, V.~Ciobanu, et~al.
\newblock {PD-1} blockade induces responses by inhibiting adaptive immune
  resistance.
\newblock \emph{Nature}, 515\penalty0 (7528):\penalty0 568--571, 2014.

\bibitem[Van~Allen et~al.(2015)Van~Allen, Miao, Schilling, Shukla, Blank,
  Zimmer, Sucker, Hillen, Geukes~Foppen, Goldinger, et~al.]{van2015genomic}
E.~M. Van~Allen, D.~Miao, B.~Schilling, S.~A. Shukla, C.~Blank, L.~Zimmer,
  A.~Sucker, U.~Hillen, M.~H. Geukes~Foppen, S.~M. Goldinger, et~al.
\newblock Genomic correlates of response to {CTLA-4} blockade in metastatic
  melanoma.
\newblock \emph{Science}, 350\penalty0 (6257):\penalty0 207--211, 2015.

\bibitem[van~der Woude et~al.(2017)van~der Woude, Gorris, Halilovic, Figdor,
  and de~Vries]{van2017migrating}
L.~L. van~der Woude, M.~A. Gorris, A.~Halilovic, C.~G. Figdor, and I.~J.~M.
  de~Vries.
\newblock Migrating into the tumor: a roadmap for t cells.
\newblock \emph{Trends Cancer}, 3\penalty0 (11):\penalty0 797--808, 2017.

\bibitem[Wang and Hillen(2007)]{wang2007classical}
Z.~Wang and T.~Hillen.
\newblock Classical solutions and pattern formation for a volume filling
  chemotaxis model.
\newblock \emph{Chaos}, 17\penalty0 (3):\penalty0 037108, 2007.

\bibitem[Wieland et~al.(2018)Wieland, Kamphorst, Adsay, Masor, Sarmiento,
  Nasti, Darko, Douek, Xue, Curran, et~al.]{wieland2018t}
A.~Wieland, A.~O. Kamphorst, N.~V. Adsay, J.~J. Masor, J.~Sarmiento, T.~H.
  Nasti, S.~Darko, D.~C. Douek, Y.~Xue, W.~J. Curran, et~al.
\newblock T cell receptor sequencing of activated {CD8 T} cells in the blood
  identifies tumor-infiltrating clones that expand after {PD-1} therapy and
  radiation in a melanoma patient.
\newblock \emph{Cancer Immunol. Immunother.}, 67\penalty0 (11):\penalty0
  1767--1776, 2018.

\bibitem[Wilkie(2013)]{Wilkie2013}
K.~P. Wilkie.
\newblock A review of mathematical models of cancer--immune interactions in the
  context of tumor dormancy.
\newblock In H.~Enderling, N.~Almog, and L.~Hlatky, editors, \emph{Systems
  Biology of Tumor Dormancy}, pages 201--234. Springer New York, New York, NY,
  2013.
\newblock ISBN 978-1-4614-1445-2.
\newblock \doi{10.1007/978-1-4614-1445-2_10}.
\newblock URL \url{https://doi.org/10.1007/978-1-4614-1445-2_10}.

\bibitem[Wolchok et~al.(2017)Wolchok, Chiarion-Sileni, Gonzalez, Rutkowski,
  Grob, Cowey, Lao, Wagstaff, Schadendorf, Ferrucci,
  et~al.]{wolchok2017overall}
J.~D. Wolchok, V.~Chiarion-Sileni, R.~Gonzalez, P.~Rutkowski, J.-J. Grob, C.~L.
  Cowey, C.~D. Lao, J.~Wagstaff, D.~Schadendorf, P.~F. Ferrucci, et~al.
\newblock Overall survival with combined nivolumab and ipilimumab in advanced
  melanoma.
\newblock \emph{N. Engl. J. Med.}, 377\penalty0 (14):\penalty0 1345--1356,
  2017.

\end{thebibliography}
%% if required, the content of .bbl file can be included here once bbl is generated
%%\input sn-article.bbl

%% Default %%
%%\input sn-sample-bib.tex%

\end{document}